\documentclass[final,3p,times]{elsarticle}
\usepackage{mathbbol}
\usepackage[ruled,vlined]{algorithm2e}
\usepackage{amsmath,bm,accents}
\usepackage{amssymb}
\DeclareSymbolFontAlphabet{\amsmathbb}{AMSb}
\usepackage{mathtools}
\usepackage{mathrsfs}
\usepackage{array}% http://ctan.org/pkg/array
\usepackage{textcomp}
\usepackage{multirow}
\usepackage{amsthm}
\usepackage{color}
\usepackage{graphicx}
\usepackage{epsfig}
\usepackage{amsfonts}
\usepackage{mathrsfs}
\newcommand{\jump}[1]{\ensuremath{[\![#1]\!]}}
\usepackage{placeins}
\usepackage{caption}
\usepackage{todonotes}

\newcommand{\R}{{\rm I\!R}}

\newcommand{\Div} {\mathrm{\bm{div}\,}}

\SetKwInput{KwInput}{Input}                % Set the Input
\SetKwInput{KwOutput}{Output}

\biboptions{numbers, sort&compress}
\usepackage[many]{tcolorbox}

\newtcolorbox{squarebrackets}[2][]{%
  empty,
  breakable,
  boxsep=0pt,boxrule=0pt,left=5mm,right=5mm,
  toptitle=3mm,bottomtitle=1mm,
  title=\colorbox{white}{#2},
  coltitle=black,
  fonttitle=\bfseries,
  underlay={%
    \draw[gray!50,line width=1mm]
      ([xshift=5mm,yshift=-0.5mm]frame.north west)--
      ([xshift=0.5mm,yshift=-0.5mm]frame.north west)--
      ([xshift=0.5mm,yshift=0.5mm]frame.south west)--
      ([xshift=5mm,yshift=0.5mm]frame.south west)
      ([xshift=-5mm,yshift=-0.5mm]frame.north east)--
      ([xshift=-0.5mm,yshift=-0.5mm]frame.north east)--
      ([xshift=-0.5mm,yshift=0.5mm]frame.south east)--
      ([xshift=-5mm,yshift=0.5mm]frame.south east)
      ;
  },
  parbox=false,#1
}

\newdefinition{rmk}{Remark}
\newtheorem*{rmk*}{Remarks}
\newproof{pf}{Proof}
\newproof{pot}{Proof of Theorem \ref{thm2}}

\journal{}

\journal{Computer Methods in Applied Mechanics and Engineering}
\journal{arXiv}
\begin{document}

\begin{frontmatter}

\title{A p-adaptive, implicit-explicit mixed finite element method for reaction-diffusion problems}

%% Group authors per affiliation:
\author[local]{Mebratu Wakeni}
\author[local]{Ankush Aggarwal}
\author[local]{Lukasz Kaczmarczyk}
\author[local]{Andrew McBride}
\author[local]{Ignatios Athanasiadis}
\author[local]{Chris Pearce}
\author[local,abroad]{Paul Steinmann}
\address[local]{Glasgow Computational Engineering Centre, University of Glasgow, Glasgow, G12 8QQ, United Kingdom}
\address[abroad]{Institute of Applied Mechanics, Friedrich-Alexander University of Erlangen-Nuremberg, Paul-Gordan-Str. 3,
D-91052, Erlangen, Germany}
%\fntext[myfootnote]{Since 1880.}

%% or include affiliations in footnotes:
%\author[mymainaddress,mysecondaryaddress]{Elsevier Inc}
%\ead[url]{www.elsevier.com}

%\author[mysecondaryaddress]{Global Customer Service\corref{mycorrespondingauthor}}
%\cortext[mycorrespondingauthor]{Corresponding author}
%\ead{support@elsevier.com}

%\address[mymainaddress]{1600 John F Kennedy Boulevard, Philadelphia}
%\address[mysecondaryaddress]{360 Park Avenue South, New York}

\begin{abstract}
A new class of implicit-explicit (IMEX) methods combined with a p-adaptive mixed finite element formulation is proposed to simulate the
diffusion of reacting species. Hierarchical polynomial functions are used to construct an $H(\Div)$-conforming base
for the flux vectors, and a non-conforming $L^2$ base for the mass concentration of the species.
The mixed formulation captures the distinct nonlinearities associated with the constitutive flux equations and the reaction terms. 
The IMEX method conveniently treats these two sources of nonlinearity implicitly and explicitly, respectively,
within a single time-stepping framework. The combination of the p-adaptive mixed formulation and the IMEX method delivers 
a robust and efficient algorithm. The proposed methods eliminate the coupled effect of
mesh size and time step on the algorithmic stability. A residual based a posteriori error estimate that provides an upper bound of the 
natural error norm is derived. The availability of such estimate which can be obtained with minimal computational effort and the 
hierarchical construction of the finite element spaces allow for the formulation of an efficient p-adaptive algorithm. A series of 
numerical examples demonstrate the performance of the approach. It is shown that the method with the p-adaptive strategy accurately 
solves problems involving travelling waves, and those with discontinuities and singularities. The flexibility of the formulation
is also illustrated via selected applications in pattern formation and electrophysiology.
% of mesh size $h$ and time step $\Delta t$
\end{abstract}


\end{frontmatter}

% \newpage
% \section*{Highlights}
% \begin{itemize}
% \item A p-adaptive, mixed formulation with IMEX methods for reaction-diffusion problems is presented.
% \item The improved performance of the proposed method over a standard, single-field formulation for selected problems is demonstrated.
% \item A strong case for why the proposed mixed methods are unambiguously essential for solutions with lower
% regularities is made.
% \item The capabilities of the method using selected problems in pattern formation, and electrophysiology are demonstrated.
% \end{itemize}
% \newpage
% %\linenumbers

\section{Introduction}
\subsection{Motivation}
The spatio-temporal dynamics of multiple  species interacting through a combination of two
distinct mechanisms, namely reaction and diffusion, can be described by {\em reaction-diffusion} equations.
Reaction refers to the inter/intra species interactions, resulting in the production and  extinction of species.
It is embodied in a term that is referred to as {\em reaction kinetics} $f$, a function of the mass concentration(s)
$m$ of the involved species. Diffusion refers to the flow of substance (concentration) in space,
and it is mathematically described by a flux $\bm{h}$ related to $m$ (and/or its spatial gradient) through a
constitutive equation. Reaction-diffusion models are relevant in various important
applications, including tissue morphogenesis and pattern formation \cite{Chaplain2001387, GARIKIPATI2017192,
TAPASWI1986213, gilbert2000develop}, tissue remodelling \cite{AMBROSI2011863, Morishita2008, Tewary4298,
MEIER2008481, Ryser2010, SIMPSON2006282}, electrophysiology \cite{Kerckhoffs2006, Rubin1996}, and epidemiology
\cite{Zhang2013, Wang2011, Wilson1997}.

% While the equations in these applications look slightly different, their fundamental mathematical structure
% is that of reaction-diffusion type.

The aforementioned applications motivate the need for robust and efficient numerical methods for solving
reaction-diffusion  problems. Various numerical methods have been proposed for approximating the solutions of
reaction-diffusion problems. Meshless methods in conjunction with operator-splitting techniques were used in
\cite{HEMAMI20193644} in one- and two-dimensions. In \cite{OLMOS20092258}, Fitzhugh-Nagumo type models are solved
using a multidomain algorithm based on a pseudospectral approach.
A review of some finite difference based methods in one-dimension can be found in \cite{RAMOS1983538}. Furthermore, a
finite difference scheme was constructed for the simulation of waves in excitable media using a two-variable
reaction-diffusion equation in \cite{BARKLEY199161}. However, finite difference and spectral algorithms are suitable
only for approximations over relatively simple domains.

\subsection{Spatial discretisation}
The finite element method, due to its capabilities in handling arbitrary geometry and nonlinearities and its
strong theoretical foundation, is a natural choice for solving reaction-diffusion problems. The majority of the finite
element numerical approaches used in the literature are based on the standard, single-field
formulation. The standard formulation for reaction-diffusion was employed in a computational framework
for the coupling of reaction-diffusion and elasticity in \cite{Ricardo2015Ricardo}. In \cite{LANG1998105}, a multilevel
finite element approach with spatial and temporal adaptivity was constructed for reaction-diffusion problems.
In \cite{TUNCER201545}, a projected finite element approach on stationary closed surface geometries,
together with a backward Euler time integration, was used to simulate pattern-formation in biological applications.
In \cite{MACDONALD2016207}, a moving mesh finite element method was constructed for simulating chemotaxis
in two-dimensions. A multigrid finite element method on stationary and evolving surfaces was proposed in
\cite{Landsberg2010}.  A semi-linear multistep finite element was constructed in \cite{Mergia2020}
for the two-dimensional simulation of pattern formation in ecological application. In such standard formulations,
only the mass concentration $m$ is solved for.  The other physically important quantity, the flux $\bm{h}$, is
obtained as a post-processing step, decreasing the accuracy of its approximation.

Solutions of reaction-diffusion problems exhibit a variety of phenomena from the formation of
travelling waves to complex structures like dissipative solitons. Some solutions may even involve low regularity
features such as evolving jump discontinuities and singularities. The $H^1$-conforming basis functions used
in the standard formulation impose an unnecessarily high regularity requirement; for solutions
displaying low regularity features the approximate solution may never converges to the true solution.
Mixed finite element methods \citep[see e.g.][]{Boffi2013, FRANCA198889, ARNOLD1990281} offer an elegant solution for such problems. 
Mixed finite element methods are two-field formulations, which employ an $H(\Div)$-conforming basis for
the flux and a $L^2$-conforming basis for the mass concentration. This combination of basis functions relaxes
the conformity requirements, allowing a wider class of solutions to be approximated accurately.

% One drawback of
% the mixed formulation is that it introduces additional degrees-of-freedom. However, by exploiting the fact that
% matrices associated to the $L^2$ basis functions can be inverted in an element-by-element manner, an efficient
% block solver with exact Schur complement preconditioner can be constructed.

Numerical studies of reaction-diffusion type problems using mixed methods are, by comparison with standard formulations, relatively few.
In \cite{FU2016102}, a stabilised mixed formulation in combination with a first-order implicit time integration was proposed
for solving steady and unsteady state reaction-diffusion problems. In this approach, $H^1$-conforming finite
element spaces are used for $m$, and $L^2$-conforming spaces for $\bm{h}$. With regard to the finite element space
used for $m$, such a method has no particular advantage over the standard finite element method in terms of accuracy.
A two-grid approach based on a variation of a mixed method with an implicit temporal integration was proposed
and analysed in \cite{liu_chen_2017}.
\subsection{Temporal discretisation}
Most numerical procedures for reaction-diffusion equations that utilise finite elements, usually approach the temporal integration
using either fully-implicit or fully-explicit methods. It is well-established that explicit methods can be very
efficient and are easy to implement, however, they usually suffer in terms of algorithmic stability, and impose severe time
step restriction arising from the diffusion term \cite{Ruuth1995}. Implicit methods are known for
their greater stability, but can be challenging in terms of implementation, and are usually less efficient as they lead
to the solution of a large system of algebraic equation. In addition, for nonlinear problems, it is necessary to derive
and compute tangent matrices that includes implicit nonlinearities at each time step, adding further
inefficiencies. Implicit-explicit (IMEX) methods mitigate such problems by
combining the advantages of explicit and implicit methods \cite{Ascher1995}. By treating the non-local diffusion
term (involving a spatial derivative) implicitly, and the local reaction term (without a spatial derivative)
explicitly, one can eliminate the coupling effect of spatial mesh size $h$ and the time step size $\Delta t$ on the
stability condition. This allows the spatial mesh to be refined adaptively without the need for reducing the time-step size.
The application of IMEX methods appears conducive for reaction-diffusion problems, however, most of the literature
on IMEX methods for such problems are limited to classical spatial discretisation techniques, such as finite-difference
\cite{Zhang2015, Ruuth1995, Farago2013} and standard finite element \cite{Lakkis2013,Mergia2020,LIN2020124944}.

The mixed method allows for the nonlinearities that may appear in the flux constitutive equation and the reaction term
to be considered separately.  For stability reasons, the flux constitutive equation must
be treated implicitly, while the nonlinearities in the reaction term can be handled explicitly. In this presentation we present a robust and efficient numerical algorithm based on mixed formulation with IMEX temporal integration methods
for problems of reaction-diffusion type.

% The implementation of the proposed algorithm is performed on the open-source library MoFEM \cite{Kaczmarczyk2020}, utilising finite element basis
% functions based on hierarchical construction \cite{ainsworth2003hierarchic}. MoFEM in turn bases its functionalities on two other open-source libraries MOAB,
% mesh-oriented database, for handling mesh related data \cite{Tautges2004MOAB, Tautges2004}, and PETSc for managing parallelisation
% associated to algebraic operations and time discretisation \cite{Balay1997, Balay2019, abhyankar2018petsc}.
\subsection{Manuscript organisation}
The contribution is organised as follows. In Section~\ref{sec:model}, a general mathematical model of multi-species
reaction-diffusion systems is presented briefly. The weak formulation of the model using a mixed approach,
is described in Section~\ref{sec:mixed-form}. In Section~\ref{sec:num-method}, relevant aspects of the numerical
procedure for the temporal discretisation using the IMEX method and the spatial approximation using mixed
Galerkin approaches are presented. Finally, in Section~\ref{sec:num-exam} the performance and capabilities of the proposed
formulations are demonstrated using various numerical examples. Here, the
performance of the mixed and the standard formulations are compared, and finally some selected examples relevant to pattern formation,
ecology, and electrophysiology are simulated using the mixed method.

\section{Model overview}\label{sec:model}
Consider $n$ species, each with mass concentration $m_i$, where $i = 1,2,3,\dots,n$, interacting in an open, bounded region
$\Omega \subset \R^{d}$ ($d = 1, 2,\textnormal{or } 3$). The local form of the mass balance, for each of the species,
is given by
\begin{equation}\label{eq:mass-balance}
 \dot{m_i}+ \Div \,\bm{h}_i = f_i(m_1,\dots,m_n),~~~~~~~~~~~ i = 1, \dots, n,
\end{equation}
\begin{figure}
\begin{center}
  \begin{tabular}{c}
  \includegraphics[angle=0,width=7.8cm,height=7.8cm ]{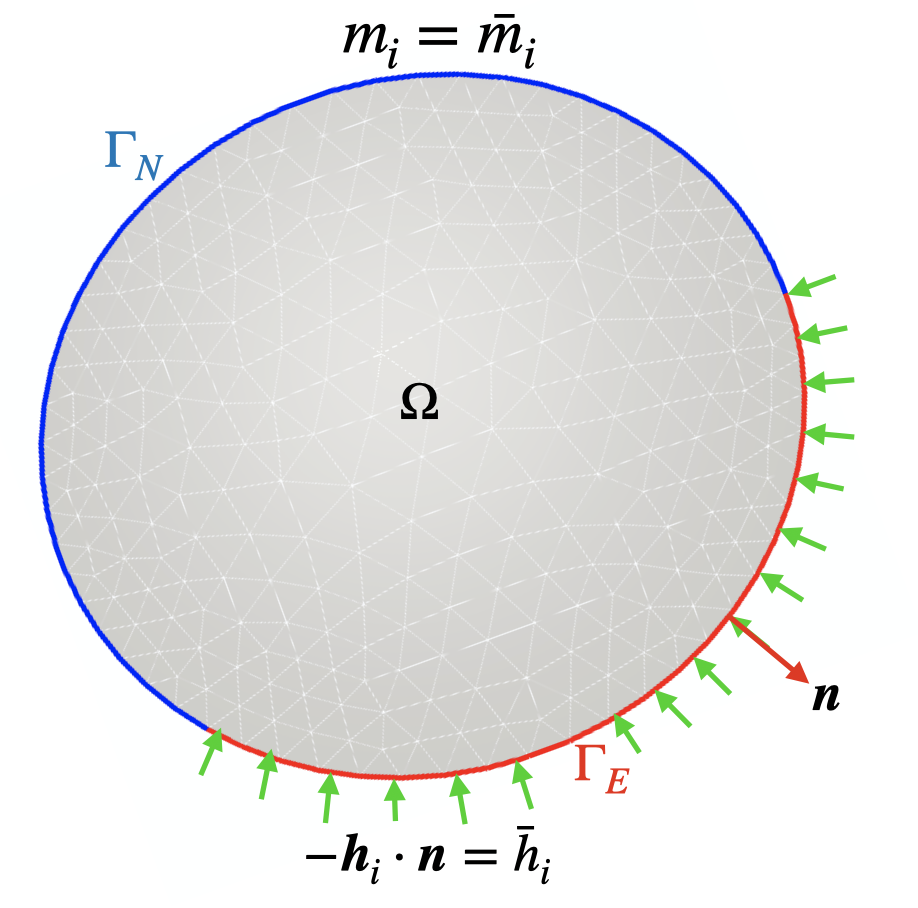}
  \end{tabular}
  \caption{Schematics of the domain $\Omega$ with its boundary partitions $\Gamma_E$ and $\Gamma_N$}
  \label{fig:scketch}
\end{center}
\end{figure}
where $\bm{h_i}$ denotes the concentration flux of the $i^{\text{th}}$ species, and $f_i$ is the chemical kinetics term that
represents  the rate of production or degradation of species concentration of the  $i^{\text{th}}$ species as a result of
its interaction with other species. In addition to the mass balance equation \eqref{eq:mass-balance}, a constitutive relation
relating the flux $\bm{h}_i$ to the mass concentration $m_i$ is required. A commonly used constitutive relation
is given by
\begin{equation}\label{eq:flux-const}
\bm{h}_i = -\bm{D}_i\nabla\, m_i,~~~~~~~~~~~ i = 1, \dots, n.
\end{equation}
Here $\bm{D}_i$ is a symmetric and positive-definite second-order tensor representing a potentially spatially varying
diffusivity/mobility of the $i^{\text{th}}$ species on the domain $\Omega$.  Let $\Gamma_N$ and
$\Gamma_E$ be nonoverlapping portions of the boundary of $\Omega$, denoted by $\Gamma$ (see Fig.~\ref{fig:scketch}), such that
$\overline{\Gamma_N \cup \Gamma_E} = \Gamma$. 
The prescribed boundary conditions imposed on these partitions are
\begin{align}
m_i &= \bar{m}_i~~~~~\text{on}~\Gamma_N\label{eq:bd-N},~\text{and}\\
-\bm{h}_i\cdot \bm{n} &= \bar{h}_i~~~~~~\text{on}~\Gamma_E \label{eq:bd-E},
\end{align}
where $\bm{n}$ represents the unit outward normal vector to the boundary $\Gamma$. A complete description of the problem
also requires the prescription of initial conditions for each $m_i$, which read as
\[
m_i(\bm{x}) = m_i^0(\bm{x}),~~~~~~\text{at}~t=0,~~\forall \bm{x} \in \Omega.
\]

\section{Weak formulations} \label{sec:mixed-form}
The focus here is on the mixed formulation. However, for the sake of completeness, the standard, single-field formulation is
first briefly stated. Thereafter, a detailed presentation of the mixed formulation and its spatial and temporal discretisation
is given.

In the context of initial-boundary value problems, such as analysed here, it is helpful to view functions of space and
time as mappings from the time interval of interest $\mathbb{I} = [0,~T]$ to the corresponding functional space.
For example, a function $u \in L^2(\Omega; \mathbb{I})$ is understood as the map $u:\mathbb{I} \to L^2(\Omega)$
(the space $L^2(\Omega)$ denotes the space of measurable functions which are square integrable over the domain $\Omega$).
In addition to the functional space $L^2(\Omega;\mathbb{I})$  we also make use of the space $H^1(\Omega; \mathbb{I})$ and
$H(\Div,\Omega; \mathbb{I})$, where
\begin{align*}
  &H^1(\Omega) = \lbrace m \in L^2(\Omega) : ~~\nabla\;m \in [L^2(\Omega)]^d\rbrace ,~~\text{and}\\
  &H(\Div,\Omega) = \lbrace \bm{h} \in [L^2(\Omega)]^d : ~~\Div\;\bm{h} \in L^2(\Omega)\rbrace .
\end{align*}
The natural norms endowed by $H^1(\Omega)$ and $H(\Div,\Omega)$ are, respectively, given by
\begin{align*}
\Vert m \Vert^2_{1,\Omega}~~ &\coloneqq   \Vert m\Vert^2_{0,{\Omega}} + \Vert \nabla\;m\Vert^2_{0,{\Omega}},~\text{and}\\
\Vert \bm{h} \Vert^2_{\Div,\Omega} &\coloneqq   \Vert\bm{h}\Vert^2_{0,{\Omega}}+ \Vert \Div\;\bm{h}\Vert^2_{0,{\Omega}},
\end{align*}
where $\Vert \cdot \Vert_{0,\Omega}$ denotes the standard $L^2$-norm for scalar or vector-valued functions.

The standard weak problem is defined as:
\begin{squarebrackets}{Standard formulation}
  find $m_i \in H^1(\Omega)$, satisfying the boundary conditions
  \eqref{eq:bd-N}, such that
  \begin{equation}\label{eq:standard}
  \dfrac{\mathrm{d}}{\mathrm{d}t}(v,~m_i)_{\Omega}+(\nabla v,~\bm{D}_i \nabla m_i)_{\Omega} = (v,~\bar{h}_i)_{\Gamma_E} +
  \ell_i(v),
  ~~~~\forall v \in H^1_{0N}(\Omega),
  \end{equation} 
\end{squarebrackets}

where $H^1_{0N}(\Omega)$ is a subspace of $H^1(\Omega)$ that contains functions whose trace on $\Gamma_N$ vanish. Here, it should be noted
that the test function $v$ is time-independent.
The functional
$\ell_i: H^1(\Omega)\to \R$ is defined by
\begin{equation}\label{eq:rhs-fun}
\ell_i(v) \coloneqq ( f_i,~ v)_{\Omega}.
\end{equation}
The pairings $(\cdot,~\cdot)_{\Omega}$ and $(\cdot,~\cdot)_{\Gamma_{N/E}}$ represent the standard $L^2$ inner product over
the domain $\Omega$ and the boundary $\Gamma_N$ or $\Gamma_E$, respectively.

As can be seen from equation \eqref{eq:standard}, the boundary condition \eqref{eq:bd-E} is incorporated
into the weak form, while equation \eqref{eq:bd-N} is enforced as a constraint on the trial solutions. Thus, for the standard
formulation, \eqref{eq:bd-E} is a natural boundary condition and \eqref{eq:bd-N} is an essential boundary condition.

For a mixed finite element formulation, in addition to $m_i$ the flux $\bm{h}_i$ is
an unknown variable. Thus, the constitutive relation \eqref{eq:flux-const} is re-written as
\begin{equation}\label{eq:flux-resist}
\bm{D}^{-1}_i\bm{h}_i + \nabla\,m_i = \bm{0}.
\end{equation}
The mixed weak form associated with equations \eqref{eq:flux-resist} and \eqref{eq:mass-balance} reads:
\begin{squarebrackets}{Mixed formulation}
find
$(\bm{h}_i,~m_i)\in H_E(\Div,\Omega;\mathbb{I}) \times L^2({\Omega;\mathbb{I}})$ such that
\begin{align}
a_i(\bm{\tau},~\bm{h}_i) - b(\bm{\tau},~m_i) &= (\bm{\tau}\cdot\bm{n},~\bar{m}_i)_{\Gamma_N},~~~~~\forall \bm{\tau}\in H_{0E}(\Div,~\Omega) \label{eq:weak-1} \\
\frac{\mathrm{d}}{\mathrm{d}t}\,c(v,~m_i) + b(\bm{h},~v) &= \ell_i(v),~~~~~~~~~~~~~~~~~~\forall v \in L^2(\Omega)\label{eq:weak-2}
\end{align}
\end{squarebrackets}
where the spaces $H_E(\Div,\Omega)$ and $H_{0E}(\Div,\Omega)$ are  subspaces of $H(\Div,\Omega)$. Also, the test functions $v$ and $\bm{tau}$ are 
time-independent.  
The vector-valued functions in the former space satisfy the boundary condition \eqref{eq:bd-E}, whereas functions in the
latter space satisfy a vanishing normal component at the boundary $\Gamma_E$.

The bilinear forms $a_i:H(\Div,\Omega)\times H(\Div,\Omega)\to \R$, $b:H(\Div,\Omega)\times  L^2(\Omega)\to\R$,
and $c:L^2(\Omega)\times L^2(\Omega)\to\R$ are defined by
\begin{align*}
a_i(\bm{\tau},~\bm{h}_i) &\coloneqq (\bm{\tau},~\bm{D}^{-1}_i \bm{h}_i)_{\Omega},\\
b(\bm{\tau},~ m_i) &\coloneqq (\Div\,\bm{\tau},~ m_i) _{\Omega},\\
c(m_i,~v) &\coloneqq (m_i,~v)_{\Omega},
\end{align*}
for any $\bm{\tau}, \bm{h}_i\in H(\Div,\Omega)$ and $m_i,v\in L^2(\Omega)$.
\begin{rmk*}
\item[1.] When the $L^2(\Omega)$ space of test functions is replaced by its discrete counterpart, the test function $v$
is chosen such that it vanishes everywhere in the domain $\Omega$ except on a given element $\Omega^e$. This in turn
implies that for $v$ in \eqref{eq:weak-2} set to unity on the element $\Omega^e$, one obtains
% that, if one chooses a test function...
\[
\frac{\mathrm{d}}{\mathrm{d}t}\int_{\Omega^e}m_i\;\mathrm{d}\Omega  +\int_{\Omega^e} \Div\,\bm{h}_i\;\mathrm{d}\Omega = \int_{\Omega^e} f_i\;\mathrm{d}\Omega.
\]
As a result, such an approximation method is said to have a locally conservative property. That is, the conservation of mass
\eqref{eq:mass-balance} is satisfied on each element.
\item[2.] The classification of boundary conditions in the mixed formulation is opposite to the standard single-field
case.  In the standard formulation, the boundary condition \eqref{eq:bd-E} is a natural one as it does not require a
priori prescription on the space of trial or test spaces. By contrast, it becomes an essential boundary condition in the
mixed formulation since the trial and test functions require the normal flux at the boundary $\Gamma_E$ to be prescribed a priori.
The role of equation \eqref{eq:bd-N} is also reversed, that is, it becomes essential in the standard formulation but natural in
the mixed formulation.
\end{rmk*}

\section{Discretisation}\label{sec:num-method}
The temporal discretisation of the weak formulations,
\eqref{eq:weak-1} and \eqref{eq:weak-2}, using a combination of implicit and explicit methods is now presented.
Then, the discrete counterparts of the
spaces $H(\Div,\Omega)$ and $L^2(\Omega)$ are detailed in the context of the hierarchical construction of
shape functions over a triangular/tetrahedral mesh.
\subsection{Temporal discretisation}
Consider first the temporal discretisation of the mixed formulation  using a class of IMEX methods.
The time interval of interest is partitioned into subintervals $[t_{n-1},~t_{n}]$
with step size $\Delta t_{n}=t_{n}-t_{n-1}$. Note that the partition need not be uniform, that is, the step-sizes need not
be equal. As equation \eqref{eq:weak-1} is without a time derivative, we treat it fully implicitly
at the current time $t_n$. For the IMEX method only the second equation \eqref{eq:weak-2} involving
time derivatives is relevant. For clarity of notation, we drop the subscript $i$ from the weak forms \eqref{eq:weak-1}
and \eqref{eq:weak-2}, and, in a general multistep context, replace each term by interpolation or extrapolation formulas
as linear combinations of previous discrete values, as defined by
\begin{align}
\frac{\mathrm{d}}{\mathrm{d}t}\,c(\bullet,~m)  &\approx  \accentset{\circ}{c}(\bullet) \coloneqq \frac{\alpha_r}{\Delta t_n}\;
c(\bullet,~m^{n}) + \sum_{j=0}^{r-1}\;\frac{\alpha_j}{\Delta t_n}\;c(\bullet,~m^{n+j-r}),  \label{eq:dot-combs}\\
b(\bm{h},~\bullet)  &  \approx \widehat{b}(\bullet) \coloneqq \beta_r\;b(\bullet,~\bm{h}^{n}) + \sum_{j=0}^{r-1}
\beta_j\;b(\bm{h}^{n+j-r},~\bullet) \label{eq:im-combs},\\
\ell(\bullet) &\approx \tilde{\ell}(\bullet) \coloneqq \sum_{j=0}^{r-1} \gamma_j \;\ell^{n+j-r}(\bullet), \label{eq:ex-combs}
\end{align}
where $\ell^{k}(\bullet)$,  $k = n-r, \dots, n-1$,
represents the family of functionals $\ell(\bullet)$ defined using the time discrete values of the reaction kinetics
$f^k = f(m^k)$, i.e., with reference to \eqref{eq:rhs-fun},
\[
\ell^k(v) = (f^k,~v)_{\Omega}.
\]
The coefficients  $\beta_0, \beta_1, \dots, \beta_r$ and $\alpha_0,\alpha_1,\dots,\alpha_r$ correspond to the implicit
interpolation formula (corresponding to equations \eqref{eq:dot-combs} and \eqref{eq:im-combs}) for the value and its time derivative of a field 
at $t_n$ on the time interval $[t_{n-r}, t_n]$.
$\gamma_0, \gamma_1,\dots,\gamma_{r-1}$ are coefficients of the explicit extrapolation (corresponding to the equation \eqref{eq:ex-combs}) of a field at $t_n$ on the
interval. The integer $r$ represents the extent to which previous step solutions, starting from the current step, are included
in the scheme. 
\begin{rmk}
Some of the commonly used IMEX schemes in the literature are: 
\begin{itemize}[]
\item IMEX schemes based on the Backward Differentiation Formula (BDF)
\begin{align*}
\text{Second-order~~~~~~~~~~~~~~~~~~~~~~~~~}& \alpha_0  = 1/2,~\alpha_1  = -2, ~\alpha_2  = 3/2,\\
                           & \beta_0 = 0, ~\beta_1 = 0, ~\beta_2 = 1,\\
                           & \gamma_0 = -1, ~\gamma_1 = 2,\\
\text{Third-order~~~~~~~~~~~~~~~~~~~~~~~~~~} & \alpha_0  = 1/24, ~\alpha_1  = -1/8, ~\alpha_2  = -7/8, ~\alpha_3  = 23/24,\\
                           & \beta_0 = 1/16, ~\beta_1 = -5/16, ~\beta_2 = 15/16, ~\beta_3 = 5/16,\\
                           & \gamma_0 = 3/8, ~\gamma_1 = -5/4, ~\gamma_2 = 15/8.        
\end{align*}
\item The second-order Crank-Nicholson -- Adams-Bashforth scheme
\begin{align*}
  & \alpha_0 = 0, ~\alpha_1  = -1, ~\alpha_2  = 1,\\
  & \beta_0 = 0, ~\beta_1 = 1/2, ~\beta_2 = 1/2,\\
  & \gamma_0 = -1/2, ~\gamma_1 = 3/2.
\end{align*}
\item The second-order additive Runge-Kutta scheme \cite{Kennedy2003, Kennedy2019}
\begin{align*}
  & \alpha_0  = -1, ~\alpha_1  = 0, ~\alpha_2  = 1\\
  & \beta_0 = 1, ~\beta_1 = 0, ~\beta_2 = 1\\
  & \gamma_0 = 0, ~\gamma_1 = 2,
\end{align*}

\end{itemize}
\end{rmk}

Substituting the discrete approximations \eqref{eq:dot-combs}-\eqref{eq:ex-combs} into the weak formulation
\eqref{eq:weak-2}, together with the discrete equation corresponding to equation \eqref{eq:weak-1} at the current
time-step renders
\begin{align}
a^n(\bm{\tau}) - b^n(\bm{\tau}) &= (\bm{\tau}\cdot\bm{n},~\bar{m})_{\Gamma_N},~~~~~~\forall \bm{\tau}\in H(\Div,~\Omega)\label{eq:bvp-1},~\text{and} \\
\accentset{\circ}{c}(v) + \;\widehat{b}(v) &= \;\tilde{\ell}(v),~~~~~~~~~~~~~~~~~~\forall v \in L^2(\Omega)\label{eq:bvp-2},
\end{align}
where
\begin{align*}
a^n(\bm{\tau}) \coloneqq a(\bm{\tau}, \bm{h}^n),~~~~ b^n(\bm{\tau}) \coloneqq  b(\bm{\tau}, ~m^n).
\end{align*}
Note, equation \eqref{eq:bvp-1} and \eqref{eq:bvp-2} constitute a boundary value problem at the
time-step $t_n$. IMEX methods can be viewed as multistep schemes involving $r-1$ previous time step's solutions. They are formally
$r$-order convergent in time.

% \begin{rmk}
% For comparison purpose, the temporal discretisation used for the standard problem
% \eqref{eq:standard} follows the IMEX discretisations \eqref{eq:dot-combs}--\eqref{eq:ex-combs}
% used for mixed discrete problems \eqref{eq:weak-1} and \eqref{eq:weak-2}, leading to the semi-discrete problems
% \eqref{eq:bvp-1} and \eqref{eq:bvp-2}.
% \end{rmk}

\subsection{Spatial discretisation}
Assume a regular decomposition $\mathcal{T}_h$ of $\Omega$ into simplexes (triangles in 2D and tetrahedral in 3D). For a
given $k\in\mathbb{Z}^+$ (a non-negative integer), denote the set of all polynomials, on a given $T\in\mathcal{T}_h$, whose
order is less than or equal to $k$ by $\mathcal{P}_k(T)$. For finite element methods which typically involve the use of
non-uniform higher-order approximations on unstructured meshes, increasing the order
of the polynomial space locally via $p$- and $hp$- adaptivity can lead to complications in enforcing global
conformity of shape functions \cite{ainsworth2003hierarchic}. Hierarchical shape functions address such problems, as well
as naturally supporting the use of $p$- and $hp$-adaptivity \cite{GUO1986, Babuska1986}. The construction
of hierarchical shape functions of arbitrary order with various conformity conditions to obtain finite element subspaces
for $L^2(\Omega)$, $H^1(\Omega)$, $H(\Div,\Omega)$ on a general unstructured meshes is detailed in Ainsworth
and Coyle \cite{ainsworth2003hierarchic}. An alternative construction of H-div conforming exact sequence element with 
arbitrary order has been proposed by Fuentes et al. \cite{FUENTES2015353}. Note, however, that space in \cite{ainsworth2003hierarchic}
consist of divergence-free zero normal functions was used in following numerical examples in Section~\ref{sec:num-exam}.

Here, hierarchic shape functions are used to define the finite element spaces corresponding to the
triangulation $\mathcal{T}_h$ such that the test spaces, for concentration $m$ and flux $\bm{h}$, are defined as
\begin{align}
\mathcal{S}_h &= \lbrace v_h\in L^2(\Omega): v_h\vert_T \in \mathcal{P}^{k}(T),\,\textnormal{where }T\in\mathcal{T}_h
\rbrace \label{eq:mass-space},~\text{and}\\
\mathcal{V}_h^{0} &=\lbrace \bm{\tau}_h\in H(\Div,\Omega): \bm{\tau}_h\vert_T \in [\mathcal{P}^{k+1}(T)]^{dim}\,
\textnormal{ and } \bm{\tau}_h \cdot \bm{n} = 0 \textnormal{ on }\Gamma_E\rbrace,\label{eq:flux-space}
\end{align}
where $dim=1,2,\text{ or }3$ refers to the spatial dimension.
While the trial space for concentration $m$ is also $\mathcal{S}_h$, the trial space
$\mathcal{V}_h$ for the flux $\bm{h}$ is given by
\begin{equation}
\mathcal{V}_h =\lbrace \bm{\tau}_h\in H(\Div,\Omega): \bm{\tau}_h\vert_T \in [\mathcal{P}^{k+1}(T)]^{dim}\,
\textnormal{ and } -\bm{\tau}_h \cdot \bm{n} = \bar{h} \textnormal{ on }\Gamma_E\rbrace.\label{eq:trial-flux}
\end{equation}
It should be noted that to obtain a stable pair $(m_h,~\bm{h}_h)$ the order of
approximation for $\mathcal{V}_h$ is required to be at least one order higher than that of $\mathcal{S}_h$, see,
for example \cite{Boffi2013}.

For the approximation of $(\bm{h},~m)$, we employ the finite dimensional trial space $(\mathcal{V}_h,~\mathcal{S}_h)$ and
test space  $(\mathcal{V}_h^0,~\mathcal{S}_h)$, as defined in equations \eqref{eq:mass-space}-\eqref{eq:trial-flux}, following a
Galerkin approach. Having specified the corresponding finite element spaces, the spatio-temporal discrete form of the 
equations \eqref{eq:bvp-1} and \eqref{eq:bvp-2} assumes a block matrix system given by
\begin{equation}\label{eq:matrix-sys}
\begin{bmatrix}
\mathbf{K} & ~~~~\mathbf{B}\\
\mathbf{B}^{\mathrm{T}} & -\sigma\mathbf{M}
\end{bmatrix}
\begin{bmatrix}
\mathbf{H}^n\\
\mathbf{m}^n
\end{bmatrix} = \begin{bmatrix}
\mathbf{F}\\
\mathbf{G}
\end{bmatrix},
\end{equation}
where
\[
\mathbf{K}_{IJ} = a(\bm{\tau}^h_I,~\bm{\tau}^h_J),~~\mathbf{B}_{IL} = -b(\bm{\tau}^h_I,~v^h_L),~\text{and}~\mathbf{M}_{KL} = c(v^h_K, v^h_L).
\]
Here, $I,~J$ denote the global indices corresponding to the numbering of the basis elements of $\mathcal{V}_h$,
while $K,~L$ correspond to that of $\mathcal{S}_h$. The right hand side $\mathbf{F}$ and $\mathbf{G}$ are
given by
\begin{align*}
\mathbf{F}_I &= (\bm{\tau}_I\cdot\bm{n},~\bar{m}_h)_{\Gamma_E},~\text{and}\\
\mathbf{G}_K&= \frac{1}{\beta_r}\; \sum_{j=0}^{r-1}[\gamma_j\;\ell^{n+j-r}(v_K) - \beta_j\;b(v_K,~\bm{h}_h^{n+j-r})],~\text{respectively}.
\end{align*}
The vector $\mathbf{H}^n$ and $\mathbf{m}^n$ are the solution vectors containing the degrees of freedom (dof) associated with the current
values of $\bm{h}^n$ and $m^n$, respectively. The coefficient $\sigma = \alpha_r/[\Delta t_n\; \beta_r]$ is a shift
coefficient of the mass matrix $\mathbf{M}$.  The matrices $\mathbf{K}$ and $\mathbf{M}$
are positive definite and symmetric. Solvability of the block system \eqref{eq:matrix-sys} also requires that $\mathbf{B}$ as a linear
map is surjective (see, for example, \cite{Boffi2013} Section 3.3). The requirement that the order of
the flux shape function should be at least one order higher than the mass concentration shape function is a sufficient condition
for the surjectivity of $\mathbf{B}$.

Formally, for sufficiently smooth solutions, the expected rate of convergence for the spatial approximation employing the
finite element spaces $\mathcal{S}_h$ and $\mathcal{V}_h$ will be of order $k+1$ in both the $L^2$ and the natural
norms.

\begin{rmk*}
\item[1.] One of the most important implications of the mixed formulation, from a computational perspective, is that
the matrix $\mathbf{M}$ can be inverted locally on an element-by-element basis, and the inverse is sparse. This is due to
the fact that there is no conformity requirement on $m_h \in \mathcal{S}_h\subset L^2(\Omega)$ over element boundaries.
\item[2.] The consequence of the above observation is that one can efficiently solve the block system using a
solver that utilises a Schur complement preconditioner. More precisely, one can exactly compute the sparse Schur
complement $\mathbf{S} = \mathbf{K} + \mathbf{B}\mathbf{M}^{-1}\mathbf{B}^{\mathrm{T}}$ in an efficient manner.
\end{rmk*}
\section{A posteriori error estimators and p-adaptivity}
Adaptive finite element methods are a fundamental numerical approach in science and engineering applications. 
The success of an adaptive algorithm relies on the availability of a good error indicator (or a posteriori error estimator)
that provides an upper bound to the true approximation error, and the complexity of its implementation and computation.  In this section, we present a residual based error estimate 
that can be computed cheaply. This estimate together with the hierarchical construction of the shape functions makes the method well suited for local p-adaptivity.
Hierarchical shape functions (and dofs) are associated with mesh entities such as vertices, edges, faces, and volumes, rather than nodes.
For example, in 2D, if the local order of two adjacent faces $F_1$ and $F_2$ sharing an edge $E$ are different, say order-$k_1$ and order-$k_2$, respectively,
then to satisfy the global $H(\Div, \Omega)$-conformity one only needs to add local shape functions of order-$\max(k_1, k_2)$ on the edge $E$.
By contrast, since there is no continuity requirement for $L^2(\Omega)$ along element interfaces, the shape functions are only associated
with faces in the 2D case.
Thus, the polynomial order of $L^2$ shape functions can be set independently in each element, thereby greatly simplifying the implementation of p-adaptivity.

To underpin an effective local p-adaptivity scheme, one requires a reliable a posteriori error estimate that provides an upper bound to the true error and,
at best, decays with the same rate as the true error as the polynomial order increases uniformly.
The energy norm on $H(\Div, \Omega)\times L^2(\Omega)$, defined by

\begin{equation}\label{eq:energy_norm}
 \Vert \bm{h}\Vert_{\Div, \Omega} + \Vert m \Vert_{0,\Omega}
\end{equation}
is the appropriate norm for measuring the magnitude of approximation errors in the mixed formulation \eqref{eq:weak-1} and \eqref{eq:weak-2}.

To develop the error estimator, and for the sake of simplicity, we shall consider a one-species mixed transient problem with homogeneous boundary condition
on $m_{h,q}$, where the additional subscript used here is to denote the order of the polynomial space, that is
\begin{align}
a(\bm{\tau}_{h, k},~\bm{h}^{n}_{h,k}) - b(\bm{\tau}_{h, k},~m^{n}_{h,k-1}) &= 0,~~~~~~~~~~~~~~~~~~~~\forall \bm{\tau}_{h, k}\in \mathcal{V}_{h, k} \label{eq:fe1} \\
\sigma\,c(v_{h, k-1},~m^{n}_{h, k-1}) + b(\bm{h}^{n}_{h, k},~v_{h, k-1}) &= (g,\;v_{h, k-1})_{\Omega},~~~~~~~~\forall v_{h, k-1} \in \mathcal{S}_{h, k-1}\label{eq:fe2}
\end{align}
where $g = \tilde{f} - \sigma m^{n-1}_{h, k-1}$, where $\tilde{f}$ denotes an extrapolation of $f$ at
the previous time-step values of $m_h$ that is determined by the specific type of the IMEX scheme, the remaining term in $g$ is a contribution from
temporal discretisation of $\dot{m}_{h, k-1} = \sigma [m^n_{h, k-1} - m^{n-1}_{h,k}]$.
Recall that for stability the orders used in equations \eqref{eq:fe1} and \eqref{eq:fe2} for the flux and mass are $k$ and $k-1$, respectively ($k > 0$).
To simplify the notation, in the remainder of this section, the current unknowns $\bm{h}^{n}_{h, k}$ and $m^{n}_{h,k-1}$ are denoted as $\bm{h}_{k}$ and $m_{k-1}$.

\subsection{Residual and interface jump based error estimates}
The a posteriori error estimate corresponding to the energy norm \eqref{eq:energy_norm} of the error is derived from the following residual and jump
based errors. Consider first a sufficiently refined regular mesh $\mathcal{T}_h$ of $\Omega$, where $h$ is the mesh parameter.
\begin{itemize}
\item Element residual errors corresponding to the constitutive and conservation of mass equations: Let $K$ be an element in $\mathcal{T}_h$,
define
\begin{align}
\eta_{_{K, R, 1}} &:= \Vert \bm{D}\nabla m_{k-1} - \bm{h}_{k} \Vert_{0,K}, \label{eq:res_err_1}\\
\eta_{_{K, R, 2}} &:= \Vert \sigma\, m_{k-1} + \Div\bm{h}_{k} - g\Vert_{0,K}. \label{eq:res_err_2}
\end{align}
\item Inter-element interface jump error: Let $e$ be an edge shared by two adjacent elements $K^+$ and $K^-$, then the jump error is defined by
\begin{equation}\label{eq:jump-err}
 \eta_{_{e,J}} := h^{-1/2}\Vert \jump{m_{k-1}}\Vert_{0, e},
\end{equation}
where the jump operator is $\jump{m_{k-1}} = m_{k-1}^+ - m_{k-1}^-$, $m_{k-1}^+$ and $m_{k-1}^-$ are values of $m_{k-1}$ at the edge $e$ from the side of 
$K^+$ and $K^-$, respectively.  
\end{itemize}
The local error estimate over a given element $K \in \mathcal{T}_h$ is thus defined by summing the error contributions from \eqref{eq:res_err_1},
\eqref{eq:res_err_2} and \eqref{eq:jump-err},
\begin{equation}\label{eq:local_post}
\eta_{_K} := \bigg[\eta^2_{K, R, 1} + \eta^2_{K, R, 2} + \sum_{e \in \partial K} \eta^2_{e, J}\bigg]^{1/2},
\end{equation}
and the global error estimate is given by
\begin{equation}\label{eq:global_post}
\eta_{_{\mathcal{T}_h}} := \bigg[\sum_{K\in \mathcal{T}_h}\big[\eta^2_{_{K, R, 1}} + \eta^2_{_{K, R, 2}}\big] + \sum_{e \in \Gamma_h} \eta^2_{_{e, J}}\bigg]^{1/2}
\end{equation}
\subsection{Upper bound}
The global residual error \eqref{eq:global_post} provides an upper bound of the energy norm of the error. To show this, we make use of the following standard
estimates
\begin{itemize}
\item[1.] Optimality estimate: Let $(\bm{h},\;m) \in H(\Div, \Omega)\times L^2(\Omega)$ be the exact solution at the current time-step,
i.e., $t_n$, then there exists $C>0$, such that
\begin{equation}
\Vert\bm{h} - \bm{h}_{k}\Vert_{\Div, \Omega} + \Vert m - m_{k-1} \Vert_{0,\Omega} \leq
C \bigg\lbrace \inf_{\bm{\tau}_{k}\in \mathcal{V}_k} \Vert\bm{h} - \bm{\tau}_{k}\Vert_{\Div, \Omega} +
\inf_{{v}_{k-1}\in \mathcal{S}_{k-1}} \Vert m - v_{k-1} \Vert_{0,\Omega} \bigg\rbrace
\end{equation}
\item[2.] Energy norm error estimate of the finite element solution:
\begin{equation}\label{eq:opt_estimate}
\begin{aligned}
  \Vert \bm{h}_{k}\Vert_{\Div, \Omega} + \Vert m_{k-1}\Vert_{0,\Omega} \leq & C \bigg[\sup_{\substack{\bm{\tau}_{k}\in \mathcal{V}_{k}\\
  \Vert \bm{\tau}_k\Vert_{\Div, \Omega} = 1}}\big\lbrace
  a(\bm{\tau}_{k},~\bm{h}_{k}) - b(\bm{\tau}_{k},~m_{k-1})
  \big\rbrace \\
 & ~~~~~~+ \sup_{\substack{v_{k-1}\in \mathcal{S}_{k-1} \\ \Vert v_{k-1} \Vert_{0,\Omega}=1}}\big\lbrace
 \sigma\,c(v_{k-1},~m_{k-1}) + b(\bm{h}_{k},~v_{k-1})
 \rbrace,
 \bigg]
\end{aligned}
\end{equation}
for some $C > 0$.
\item[3.] \emph{Saturation assumption}: One of the most crucial ingredients towards the proof of an upper bound is the saturation assumption.
Roughly, it states that the error norm decreases uniformly as we increase the order of approximation by one. More precisely, let $(\bm{h}_{k+1}, m_{k})$
and $(\bm{h}_{k}, m_{k-1})$ be approximate solutions of \eqref{eq:fe1} and \eqref{eq:fe2}, then there is $0<\beta < 1$, such that
\begin{equation}\label{eq:saturate}
  \Vert\bm{h} - \bm{h}_{k+1}\Vert_{\Div, \Omega} + \Vert m - m_{k} \Vert_{0,\Omega} <
  \beta \big[\Vert\bm{h} - \bm{h}_{k}\Vert_{\Div, \Omega} + \Vert m - m_{k-1} \Vert_{0,\Omega}\big].
\end{equation}
One can construct a mesh $\mathcal{T}_h$ on which such a saturation estimate does not hold, however, for sufficiently refined
regular mesh it always hold true.
\end{itemize}
Having the above results for the upper bound, it is sufficient to show that the error between successive approximations is bounded from
above. That is, there is a constant $C>0$ such that
\begin{equation}\label{eq:rel-estimate}
\Vert\bm{h}_{k+1} - \bm{h}_{k}\Vert_{\Div, \Omega} + \Vert m_{k}  - m_{k-1} \Vert_{0,\Omega} \leq C \eta_{_{\mathcal{T}_h}}.
\end{equation}
To show this, we first note that
\begin{align}
a(\bm{\tau}_{k+1},~\bm{h}_{k+1} - \bm{h}_{k}) - b(\bm{\tau}_{k+1},~m_{k}-m_{k-1}) &= 0,~~~~~~~~\forall \bm{\tau}_{k+1}\in \mathcal{V}_{k+1} \label{eq:rel1} \\
\sigma\,c(v_{k},~m_{k}-m_{k-1}) + b(\bm{h}_{k+1}-\bm{h}_{k},~v_{k}) &= 0.~~~~~~~~\forall v_{k} \in \mathcal{S}_{k}\label{eq:rel2}
\end{align}
Hence, by the estimate \eqref{eq:opt_estimate}, we have, for some $C>0$,
\begin{equation}
\begin{aligned}
& \Vert\bm{h}_{k+1} - \bm{h}_{k}\Vert_{\Div, \Omega} + \Vert m_{k}  - m_{k-1} \Vert_{0,\Omega} \\ 
& \leq C \bigg[
  \sup_{\substack{\bm{\tau}_{k+1}\in \mathcal{V}_{k+1}\\
                  \Vert \bm{\tau}_{k+1}\Vert_{\Div, \Omega} = 1}}
\big\lbrace
  a(\bm{\tau}_{k+1},~\bm{h}_{k+1}-\bm{h}_{k}) - b(\bm{\tau}_{k+1},~m_{k} - m_{k-1})
  \big\rbrace \\
  & ~~~~~~~~~~~~~~~~~~~~~+ \sup_{\substack{v_{k}\in \mathcal{S}_{k}\\ \Vert v_k \Vert_{0,\Omega}=1}}\big\lbrace
  \sigma\,c(v_{k},~m_{k}-m_{k-1}) + b(\bm{h}_{k+1}-\bm{h}_{k},~v_{k})
  \big\rbrace
\bigg].
\end{aligned}
\end{equation}
Now, since
\begin{equation}
  a(\bm{\tau}_{k},~\bm{h}_{k+1}-\bm{h}_{k}) - b(\bm{\tau}_{k},~m_{k} - m_{k-1}) = 0, ~~~~ \forall \bm{\tau}_{k}\in \mathcal{V}_{k},
\end{equation}
it follows for every $\bm{\tau}_{k}$ that
\begin{align*}
&a(\bm{\tau}_{k+1},~\bm{h}_{k+1}-\bm{h}_{k}) - b(\bm{\tau}_{k+1},~m_{k} - m_{k-1}) \\
& = a(\bm{\tau}_{k+1} - \bm{\tau}_{k},~\bm{h}_{k+1}-\bm{h}_{k}) - b(\bm{\tau}_{k+1} - \bm{\tau}_{k},~m_{k} - m_{k-1})\\
& = - a(\bm{\tau}_{k+1} - \bm{\tau}_{k},~\bm{h}_{k}) + b(\bm{\tau}_{k+1} - \bm{\tau}_{k},~m_{k-1})\\
& =  \sum_{K}\bigg\lbrace (\bm{\tau}_{k+1} - \bm{\tau}_{k},~\bm{h}_{k})_{K} + (\Div(\bm{\tau}_{k+1} - \bm{\tau}_{k}),~m_{k-1})_{K}\bigg\rbrace\\
& = \sum_{K} \bigg\lbrace -(\bm{\tau}_{k+1} - \bm{\tau}_{k},~\bm{h}_{k})_{K} + (\bm{\tau}_{k+1} - \bm{\tau}_{k},~\nabla m_{k-1})_{K} \bigg\rbrace+ \sum_{e\in \Gamma_h} ([\bm{\tau}_{k+1} - \bm{\tau}_{k}]\cdot \bm{n},~\jump{m_{k-1}})_{e}\\
& = \sum_{K} (\bm{\tau}_{k+1} - \bm{\tau}_{k},~\bm{h}_{k}-\nabla m_{k-1})_{K} + \sum_{e\in \Gamma_h} ([\bm{\tau}_{k+1} - \bm{\tau}_{k}]\cdot \bm{n},~\jump{m_{k-1}})_{e}\\
& \leq \sum_{K} \Vert\bm{\tau}_{k+1} - \bm{\tau}_{k} \Vert_{0, K}\;\Vert\bm{h}_{k}-\nabla m_{k-1} \Vert_{0, K}+ \sum_{e\in \Gamma_h} \Vert[\bm{\tau}_{k+1} - \bm{\tau}_{k}]\cdot \bm{n} \Vert_{0, e}\; \Vert\jump{m_{k-1}}\Vert_{0,e}\\
& = \sum_{K} \eta_{_{K, R, 1}}\Vert\bm{\tau}_{k+1} - \bm{\tau}_{k} \Vert_{0, K} + \sum_{e\in \Gamma_h}\eta_{_{e, J}}h^{1/2}\Vert[\bm{\tau}_{k+1} - \bm{\tau}_{k}]\cdot \bm{n} \Vert_{0, e}.
\end{align*}
Hence, we obtain that
\begin{equation}
  \sup_{\substack{\bm{\tau}_{k+1}\in \mathcal{V}_{k+1} \\
  \Vert \bm{\tau}_{k+1}\Vert_{\Div, \Omega} = 1}}
\big\lbrace
a(\bm{\tau}_{k+1},~\bm{h}_{k+1}-\bm{h}_{k}) - b(\bm{\tau}_{k+1},~m_{k} - m_{k-1})
\big\rbrace \leq C_1 \bigg\lbrace\sum_{K} \eta^2_{_{K, R, 1}} + \sum_{e}\eta^2_{_{e, J}}\bigg\rbrace^{1/2}
\end{equation}
for some constant $C_1 > 0$. A similar argument also leads to
\begin{equation}
\sup_{\substack{v_{k}\in \mathcal{S}_{k} \\ \Vert \Vert_{0,\Omega}=1}}\big\lbrace
  \sigma\,c(v_{k},~m_{k}-m_{k-1}) + b(\bm{h}_{k+1}-\bm{h}_{k},~v_{k})
  \big\rbrace \leq C_2 \bigg\lbrace\sum_{K}\eta^2_{_{K, R, 2}}\bigg\rbrace^{1/2}.
\end{equation}
Therefore, for $C = \max(C_1, C_2)$ we obtain the estimate \eqref{eq:rel-estimate}. Employing the saturation estimate \eqref{eq:saturate}
and \eqref{eq:rel-estimate}, it then follows that
\begin{equation}
  \Vert\bm{h} - \bm{h}_{k}\Vert_{\Div, \Omega} + \Vert m - m_{k-1} \Vert_{0,\Omega} \leq \frac{C}{1-\beta}\eta_{_{\mathcal{T}_h}}.
\end{equation}
Here the constant $C$ depends only on the approximation order $k$ and $h$.

\section{Adaptive p-refinement strategy}
Once the problem is solved with a given distribution of polynomial orders over the mesh entities, and the local a posteriori
error $\eta_K$ over each element $K \in \mathcal{T}_h$ is calculated, the next step is to apply a p-refinement strategy inspired by the well-known bulk-chasing
D\"orfler's criterion \cite{Dorfler2007}. The refinement algorithm is characterised by two parameters $\theta_{\mathrm{min}}$ and $\theta_{\mathrm{max}}$, where
$0 \leq \theta_{\mathrm{min}} < \theta_{\mathrm{min}} \leq 1$, and is performed in three stages:
\begin{itemize}
\item[]{\bf Stage 1.} Given a posteriori error estimate $\eta_K$ on each element $K\in \mathcal{T}_h$ and $\eta_{\mathrm{MAX}}
= \max_{K\in\mathcal{T}_h}\lbrace \eta_K \rbrace$,
the polynomial order over element $K$ is raised by one if
\[
\eta_K \geq \theta_{\mathrm{max}}\eta_{\mathrm{MAX}},
\]
or reduced by one if
\[
\eta_K \leq \theta_{\mathrm{min}}\eta_{\mathrm{MAX}}.
\]
 After applying this first stage of the adaptive process it may happen that the polynomial order distribution over adjacent elements 
be greater than one order. Numerical experiments (not presented here) reviled that hetrogeneity of polynomial order distribution results in
undesirable oscillatory feature of the approximated solution. Hence, following this step, certain smoothing of polynomial order over the
mesh is required, which leads us to the next stage.

\item[]{\bf Stage 2.} To smooth the polynomial order distribution, we force the difference in polynomial order between two adjacent
elements $K$ and $K'$ to not exceed one, by resetting the order on the element with smaller degree to that of with the higher degree minus one. 
That is, suppose $\mathrm{order}(K) + 1 < \mathrm{order}(K')$, then  we reset $\mathrm{order}(K) \coloneqq \mathrm{order}(K) - 1$. 
\item[]{\bf Stage 3.} This stage is responsible to maintain the $H(\Div)$-conformity of the space of flux functions after we execute the
above two stages. For each interface entity $E$ shared by two elements $K$ and $K'$, we set the order as the maximum of the
polynomial orders over $K$ and $K'$.
\end{itemize}
The adaptive p-refinement algorithm consisting of the above three stages is summarised in Algorithm~\ref{algorithm}. Following the above 
p-adaptive stages, one also needs to adjust the quadrature rules  over the mesh entities appropriately in order to match the polynomial
order distributions optimally.

\begin{algorithm}[H]\label{algorithm}
\SetAlgoLined
\KwInput{$\eta_K$ on each $K\in \mathcal{T}_h$, , $\theta_{\mathrm{max}}$ and $\theta_{\mathrm{min}}$}
 {\bf Stage 1.} Setting order on each $K$\;
  \For{$K\in \mathcal{T}_h$}{

   \If{$\eta_K\geq \theta_{\mathrm{max}}\eta_{\mathrm{MAX}}$}{
   raise polynomial order on $K$ by one\;
   }
   \If{$\eta_K\leq \theta_{\mathrm{min}}\eta_{\mathrm{MAX}}$}{
    decrease polynomial order on $K$ by one but not less than a minimum order, say 1\;
  }
 }
 {\bf Stage 2.} Order smoothing\;
 \For{$E \in \Gamma_h$ shared by two elements $K$ and $K'$ in $\mathcal{T}_h$ such that the order in $K$ is greater than that of $K'$ by more than 1}{
  set: order in $K'$ equals order in $K$ minus 1\;
 }
 {\bf Stage 3.} Setting order on the interfaces\;
 \For{ $E \in \Gamma_h$}{
   \eIf{$E$ is on the boundary with only one adjecent element $K$}{
     set order on $E$ to be equal to that of on $K$\;
   }{
    find adjacent elements $K$ and $K'$ sharing $E$\;
    Set order on $E$ to be maximum of orders on $K$ and $K'$\;
   }
 }

\caption{Adaptive p-refinements algorithm}
\end{algorithm}

\section{Numerical examples}\label{sec:num-exam}

Two groups of numerical examples are presented: the first
compares the convergence of results of the mixed scheme and the standard single-field formulation in approximating
important aspects of the solution. These include solutions involving singularities, and computation of the speed of travelling
wave solutions. Suitability of the p-adaptive mixed formulation in terms of the features of the solution
is also investigated. The second group showcases the capabilities of the p-adaptive, IMEX mixed formulation in simulating
problems of practical importance: pattern formation and electrophysiology. We investigated several IMEX schemes, for the examples
presented in this section we opted for the second-order additive Runge-Kutta scheme.

The computer implementation of the proposed numerical scheme is carried out using the open-source library MOFEM
\cite{Kaczmarczyk2020}. The library integrates and utilises other open-source libraries such as MOAB, a mesh-oriented
database \cite{Tautges2004MOAB, Tautges2004}, and PETSc \cite{Balay1997, Balay2019}. The MOAB library is used to store and manage mesh related data,
while PETSc is used for parallel operations involving linear algebra.

The IMEX methods presented in Section \ref{sec:num-method} are implemented using the PETSc (Portable, Extensible Toolkit for Scientific
computations \cite{Balay1997, Balay2019}) time solvers \cite{abhyankar2018petsc}. 
\subsection{Convergence tests}
Two cases are considered. In the first, a
spatially smooth solution is considered with a piece-wise temporal profile that stabilises after some specified time.
Thus, the approximation error after a sufficiently long simulation time is associated entirely with the spatial discretisation.
The second case considers the approximation of a one-species Fisher's type problem on a square domain $\Omega$ with
heterogeneous diffusivity.

\subsection*{a) Smooth manufactured solution}
\begin{figure}
\begin{center}
 \begin{tabular}{c}
 \includegraphics[angle=0,width=7.5cm,height=6.0cm ]{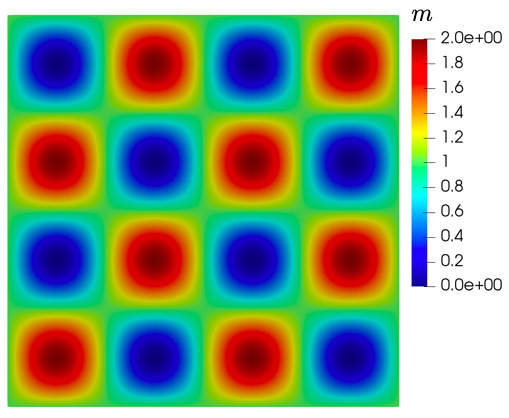}
\end{tabular}
\caption{Exact (manufactured) solution for $m$ at $t = 1$.}
\label{fig:analytic}
\end{center}
\end{figure}

It is well-known that both standard and the mixed finite element formulations are optimal in terms of convergence in the $L^2$-norm, i.e., $\mathcal{O}(h^{p+1})$,
for sufficiently smooth solutions, where $p$ is the order of the finite element space. Noting that the mesh size parameter $h$ is inversely
proportional to the number of degrees-of-freedom to the power $dim$, where $dim = 1,2\text{ or } 3$ is the space dimension, these optimality
results are confirmed practically, as shown in Fig. \ref{fig:analytic}, by considering a manufactured solution based on the smooth
function
\begin{equation}\label{eq:smooth_solution}
g(x, y) = 1 + \sin(2\pi x)\sin(2\pi y),~~~~~\textnormal{for }(x, y)\in\Omega.
\end{equation}

Consider first a one-species reaction-diffusion system over the domain $\Omega = [-1, 1]^2$ (so that $g$ vanishes on the boundary) with
isotropic diffusivity $\bm{D} = d\bm{I}$, $d = 1$, and then assume the exact (manufactured) solution for the mass concentration $m$
\begin{equation}\label{eq:manufactured}
m(x, y, t) =
\begin{cases}
  t\;g(x, y), & t< t^*, \\
 g(x, y), & t \geq t^*.
 \end{cases}
\end{equation}
for some given $t^*$. The right-hand-side source term $f$  is given by the residual of the exact solution, i.e.,
\[
f \coloneqq \dot{m} + \Div(\bm{D}\nabla m ).
\]
Note that $\dot{m}$ and $\Div(\bm{D}\nabla m )$ are also piecewise in time, i.e.,
\[
\dot{m} = \begin{cases}
g(x, y), & t< t^*, \\
0, & t > t^*,
\end{cases}~~~\textnormal{and }~
\Div(\bm{D}\nabla m) =
\begin{cases}
t\;\Div(\bm{D}\nabla g), & t<t^*,\\
~~\Div(\bm{D}\nabla g), & t\geq t^*.
\end{cases}
\]

Consequently, the source term is also temporally piece-wise which stabilises to  a time-independent profile after $t^*$.
Note also that for $t<t^*$, $\bar{m} = t$ on the boundary $\Gamma$, and for $t>t^*$, $\bar{m} = 1$.
With this set up, the temporal discretisation error after a time $t$ sufficiently greater than $t^*$ will be negligible, and the
total error is dominated by the spatial approximation. In other words, it amounts to the approximation of the steady state case
($\dot{m} = 0$) with the manufactured solution $m = g$.

\begin{figure}
\begin{center}
 \begin{tabular}{cc}
 \includegraphics[angle=0,width=8.0cm,height=6.5cm ]{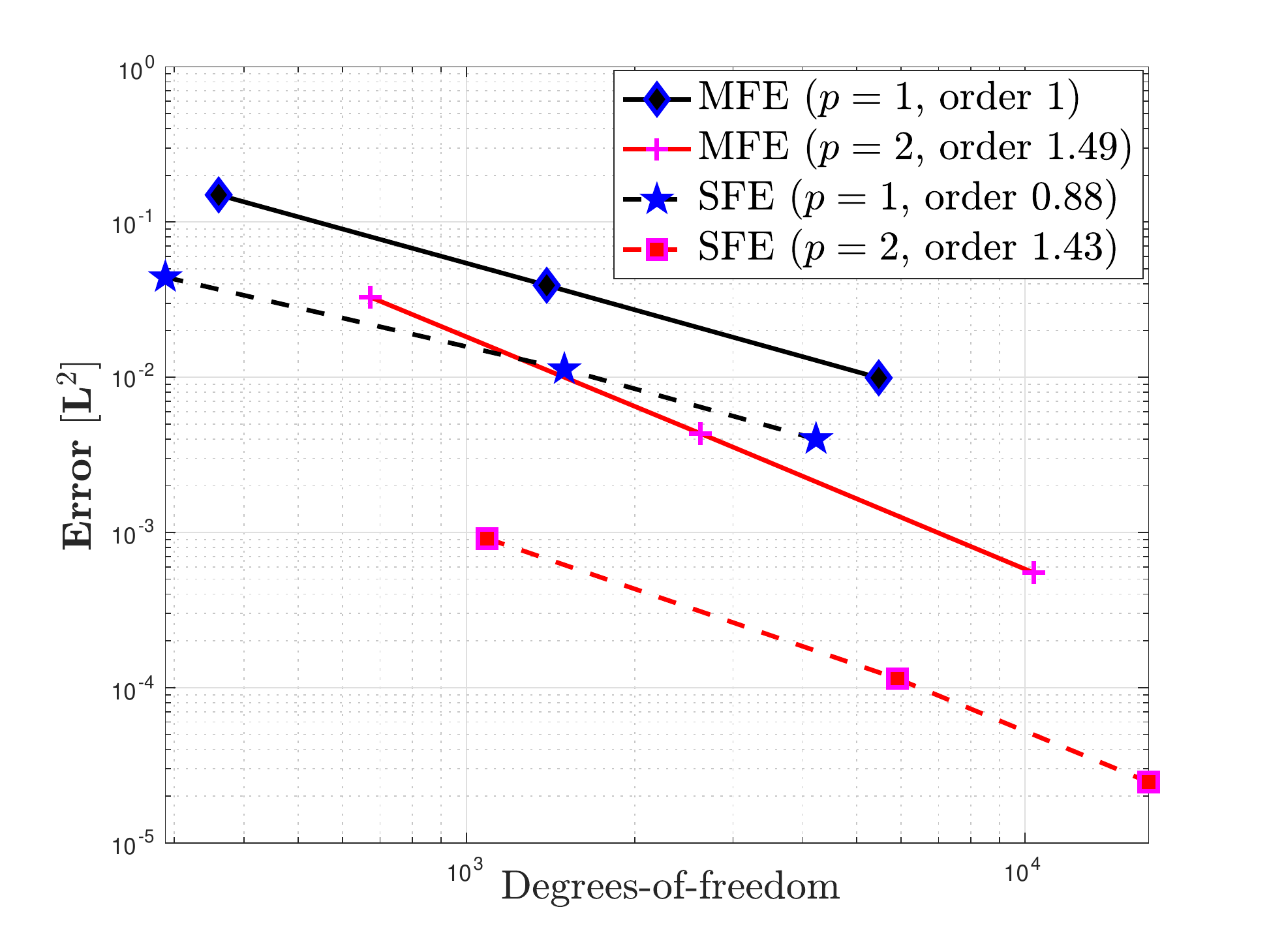} &
 \includegraphics[angle=0,width=8.0cm,height=6.5cm ]{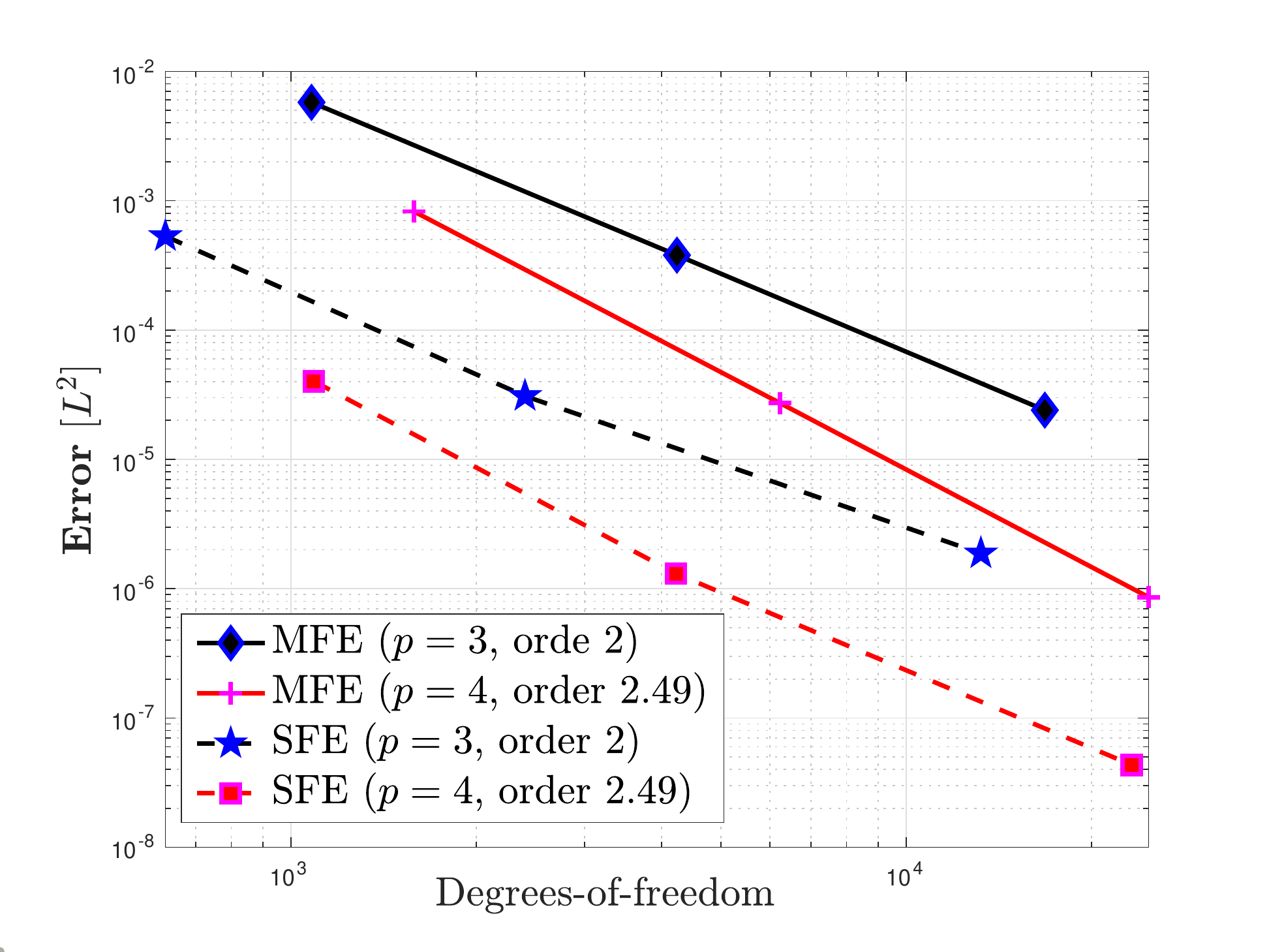}\\
 (a) & (b)
\end{tabular}
\caption{Comparison of convergence rates of the mixed (MFE) and standard (SFE) formulations with respect to the $L^2$ error norm.
Figure (a) is for orders of approximation  $p=1, 2$, while (b) is for $p=3,4$. The legend `order' stands for the absolute
value of the slope of the convergence curve once a consistent slope is established between consecutive refinements.}
\label{fig:l2-error}
\end{center}
\end{figure}

\begin{figure}
\begin{center}
 \begin{tabular}{cc}
 \includegraphics[angle=0,width=8.0cm,height=6.5cm ]{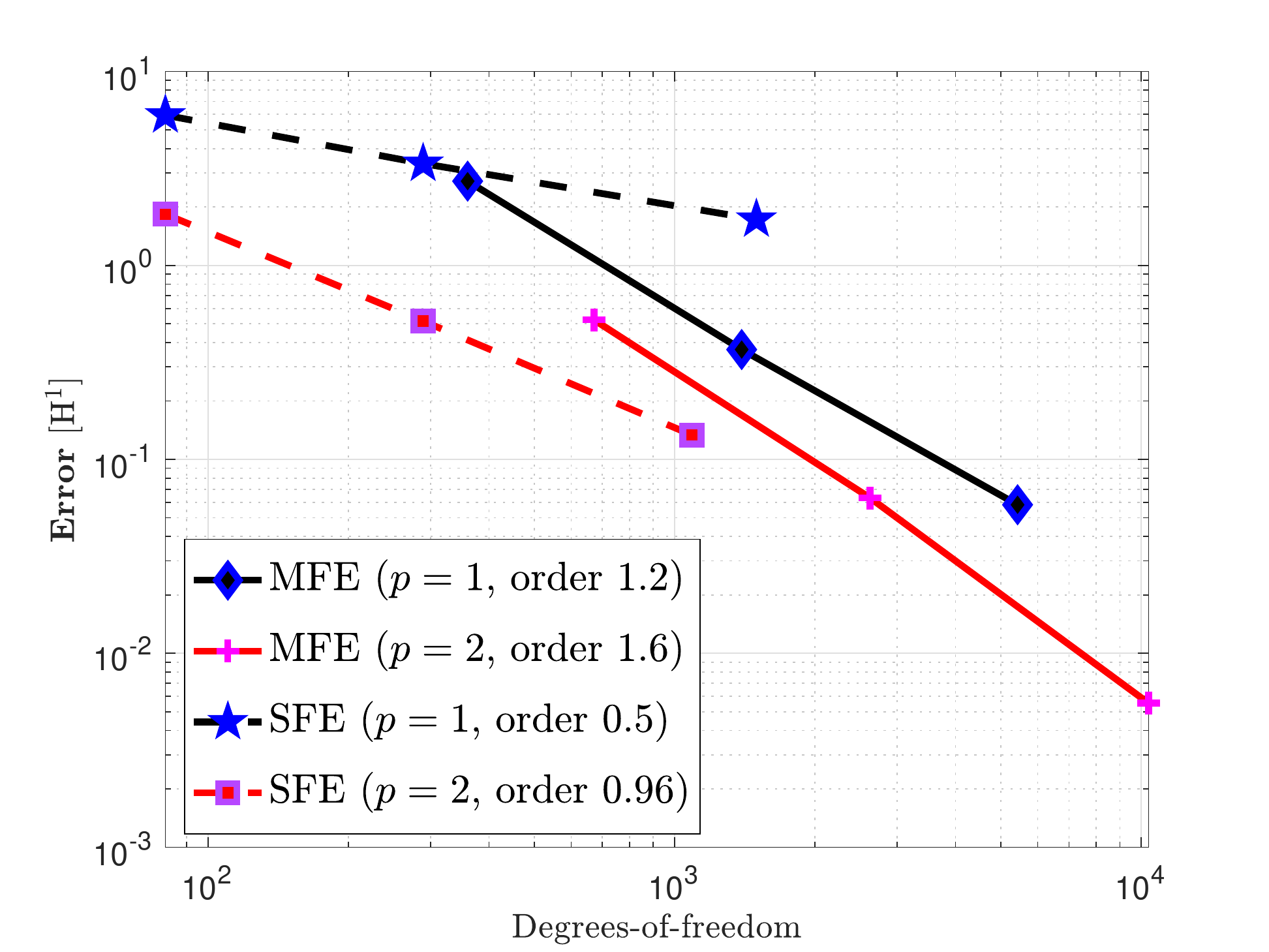} &
 \includegraphics[angle=0,width=8.0cm,height=6.5cm ]{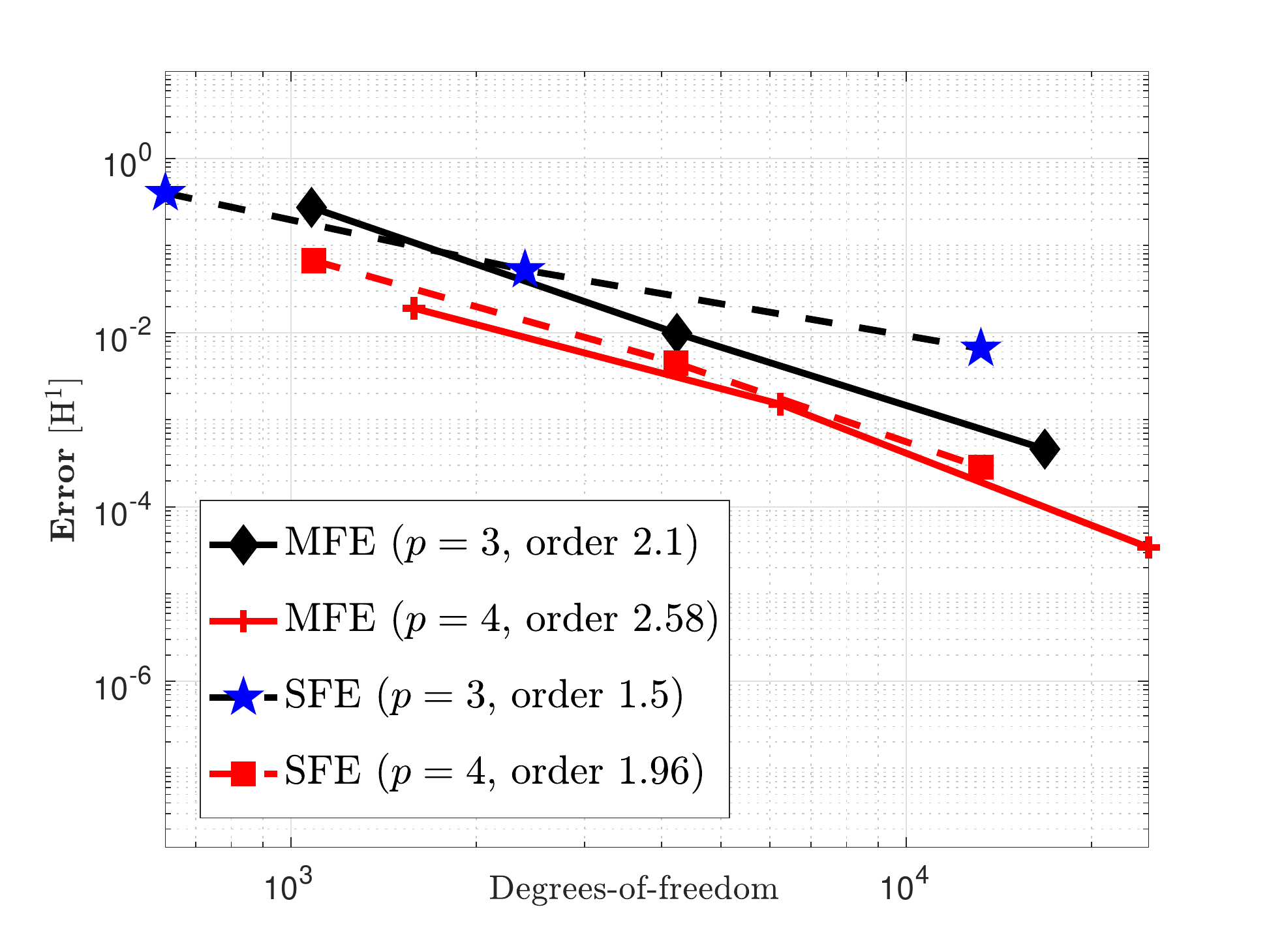}\\
 (a) & (b)
\end{tabular}
\caption{Comparison of convergence rate the mixed (MFE) and standard (SFE) formulations with respect to the $H^1$-norm.
Figure (a) is for order of approximation $p=1, 2$, for $m$ while (b) is for $p=3,4$. The legend `order' stands for
the absolute value of the slope of the convergence curve once a consistent slope is established between consecutive
refinements.}
\label{fig:H1-error}
\end{center}
\end{figure}

The convergence results presented in Fig. \ref{fig:l2-error} and \ref{fig:H1-error} are obtained by a successive refinement of
an initial uniform mesh with $h = 2/5$, and the inflection time $t^*$ and a uniform time step length $\Delta t$ are chosen
to be $1$ and $0.1$, respectively, corresponding to each mesh. The simulations are run up to $t = 10$; a sufficiently long  time to ensure that the temporal discretisation error is negligible.

With the same construction of the manufactured solution \eqref{eq:manufactured}, the convergence rate of the mixed formulation with
respect to the $H^1$-norm, given by
\[
\big[\Vert m - m^h\Vert_{0,\Omega}^2 + \Vert \bm{h} - \bm{h}^h\Vert_{0,\Omega}^2\big]^{1/2},
\]
is expected to be one order higher than that of the standard formulation, as demonstrated in Fig. \ref{fig:H1-error}.
This is due to the fact that the flux $\bm{h}_{h}$ for the standard formulation is obtained by postprocessing from $m_h$, unlike the
mixed formulation, wherein the flux is directly approximated as a primary field variable.

% Fig.~\ref{fig:l2-error} (a) and (b) compare the rate of convergence with respect to  the $L^2$-error norm of the mass
% approximations using the mixed and standard formulations. The optimal convergence rate, which is
% $\mathcal{O}(h^{k+1})$, is demonstrated by both schemes. However, the standard formulation does not directly
% approximate the flux but derives it from the approximated mass concentration. The rate of
% convergence with respect to the  $H^1$-error norm, which is given by
% \[
% \big[\Vert m - m^h\Vert_{0,\Omega}^2 + \Vert \bm{h} - \bm{h}^h\Vert_{0,\Omega}^2\big]^{1/2},
% \]
% for the mixed formulation, is expected to be one order higher than that of the standard formulation, as demonstrated
% in Fig.~\ref{fig:H1-error}.

\begin{figure}
\begin{center}
 \begin{tabular}{c}
 \includegraphics[angle=0,width=8.0cm,height=6.5cm ]{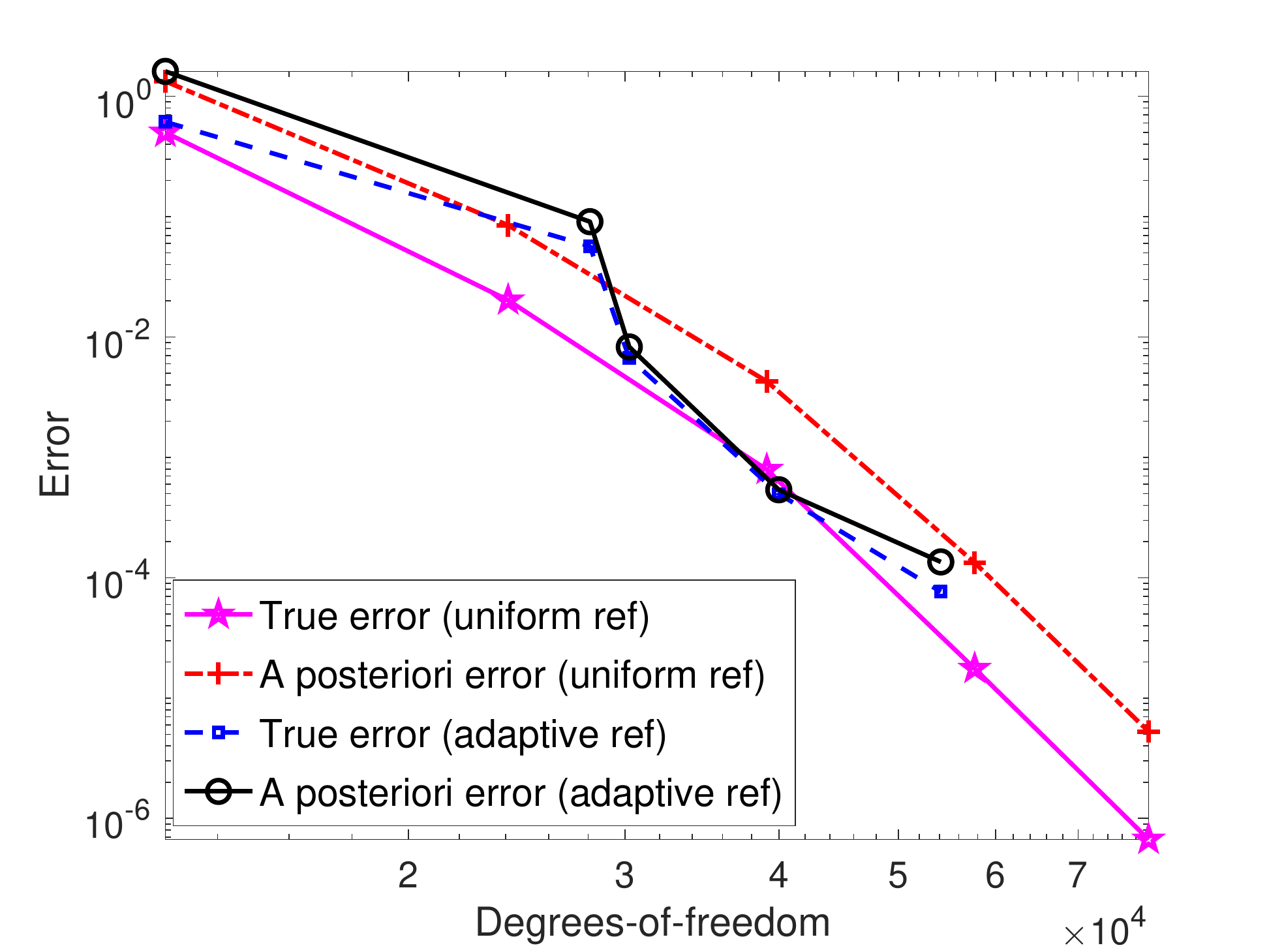}
\end{tabular}
\caption{Comparison of the p-adaptive and uniform refinement strategies for smooth and slowly varying solution.}
\label{fig:uniform_adaptive}
\end{center}
\end{figure}

\begin{figure}
\begin{center}
  \begin{tabular}{c}
  \includegraphics[angle=0,width=17.0cm,height=7cm ]{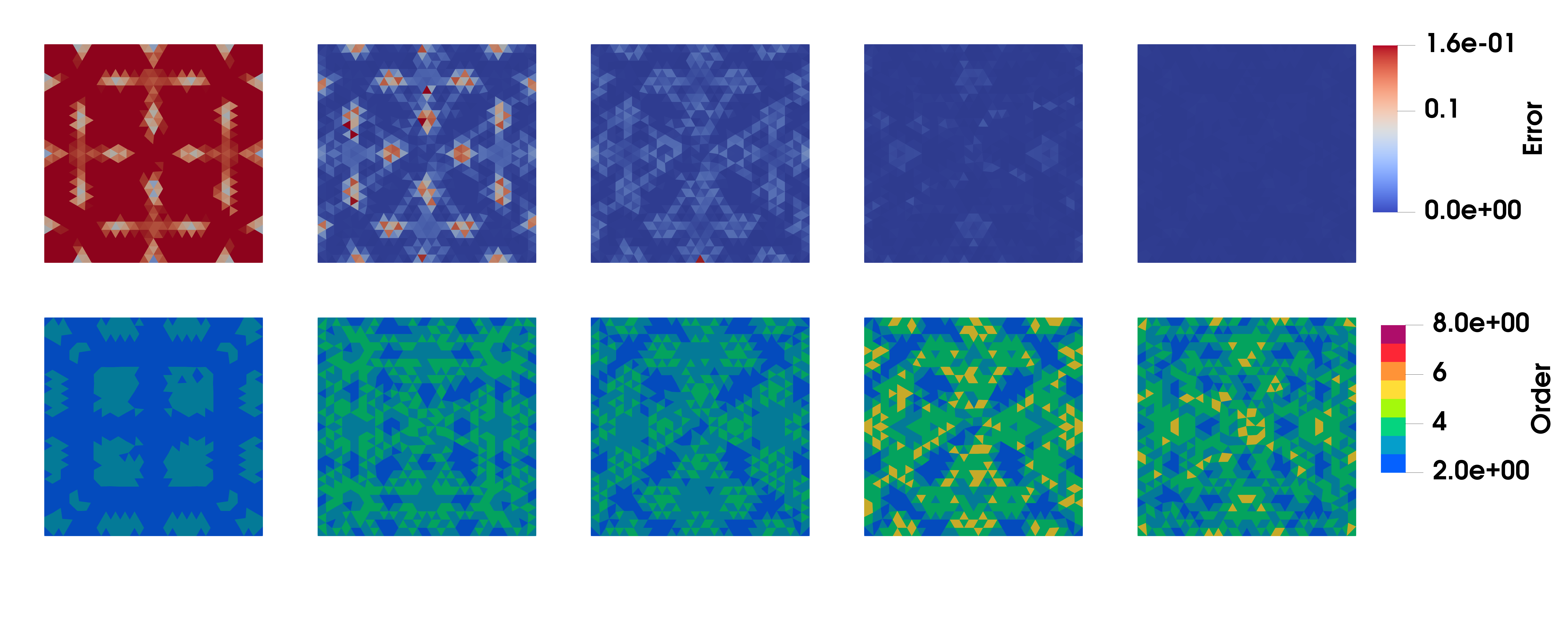}
\end{tabular}
\caption{Distribution of error and polynomial order over the mesh at each p-adaptive iteration}
\label{fig:uniform_sequence}
\end{center}
\end{figure}

The next set of examples in this group aims at demonstrating the effectiveness of the p-adaptive mixed formulation in
resolving fine features of solutions efficiently. For smooth solutions such as \eqref{eq:smooth_solution}, the variability
of the solution is almost uniform on the larger scale. In this case, the application of p-adaptivity is less effective since
the error is distributed almost uniformly. This is demonstrated in the convergence result displayed in Fig.~\ref{fig:uniform_adaptive}.
It shows that the p-adaptive strategy with parameters $\theta_{\mathrm{min}} = 0.02$ and $\theta_{\mathrm{max}} = 0.8$, representing a quite conservative
adaptive strategy, produces a convergence trend which is not generally better than that of the uniform p-refinement. As expected, at each
adaptive step, as shown in Fig.~\ref{fig:uniform_sequence}, the error is distributed almost uniformly, which leads to the marking of most of
the elements for refinement. This corresponds to the convergence result shown in Fig.~\ref{fig:uniform_sequence}, which is  not better
than the uniform p-adaptive strategy.

\begin{figure}
  \begin{center}
    \begin{tabular}{c}
    \includegraphics[angle=0,width=15.0cm,height=5cm ]{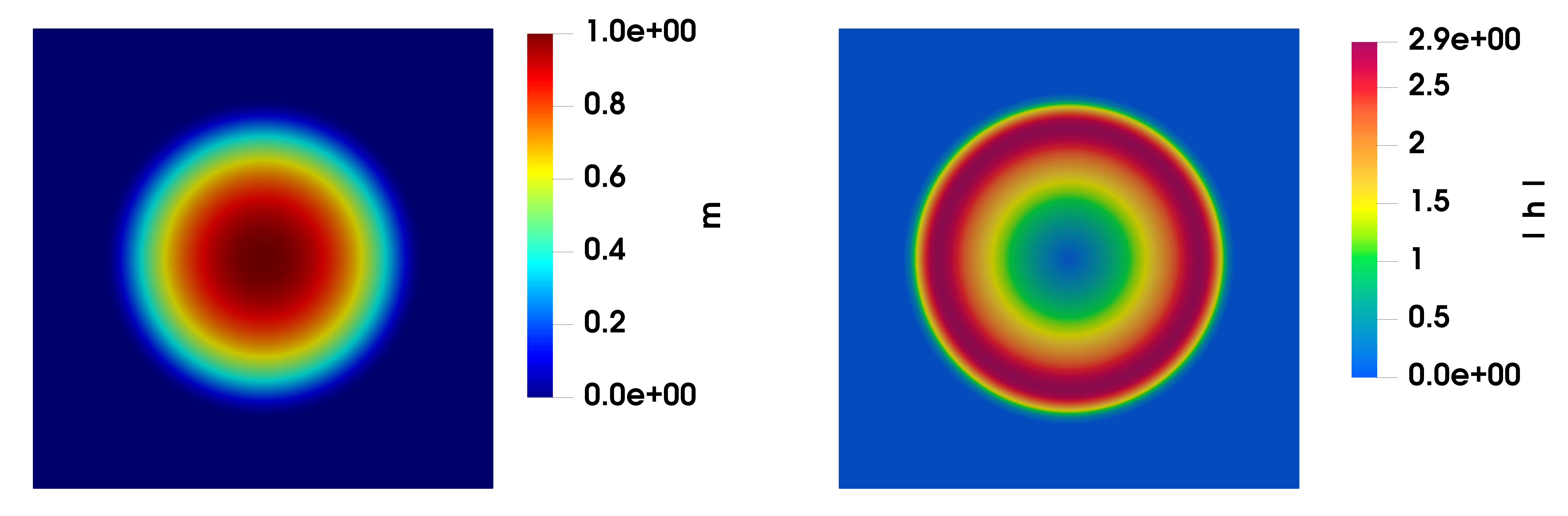}
  \end{tabular}
  \caption{Bump function (left) and norm of the gradient (right)}
  \label{fig:non_uniform}
  \end{center}
  \end{figure}

In contrast, when the solution is characterised by the presence of sudden spatial changes over the domain,
such as travelling waves, the p-adaptive algorithm becomes most effective. To demonstrate this, we consider a smooth
analytical solution (manufactured) replacing the $g$ in equation \eqref{eq:smooth_solution} by the bump function over the square domain $\Omega = [-1, 1]^2$,
\begin{equation}
g(x, y) = \begin{cases}
          \exp{\big(\frac{-x^2 - y^2}{r^2 - x^2 - y^2}\big)} & \text{if } x^2 + y^2 < r^2,\\
          0                                                 & \text{otherwise}.
          \end{cases}
\end{equation}
where $r = 0.75$. As in the previous examples, the diffusion parameter is chosen such that $d = 1$, and $t^* = 1$ and $\Delta t = 0.1$. The error
is computed at $t = 10$, as in the previous simulation, as it is far enough from the inflection time $t^* = 1$ so that the temporal discretisation
error becomes negligible. The spatial mesh is unstructured and relatively coarse. It can be easily seen that the bump function is infinitely many times differentiable hence smooth. However,
as can be seen from Fig.~\ref{fig:non_uniform}, along the circle $x^2+y^2 = r^2$ the solution changes drastically from zero
to some finite non-zero value within a relatively small distance in the radial direction.

\begin{figure}
\begin{center}
 \begin{tabular}{c}
 \includegraphics[angle=0,width=8.0cm,height=6.5cm ]{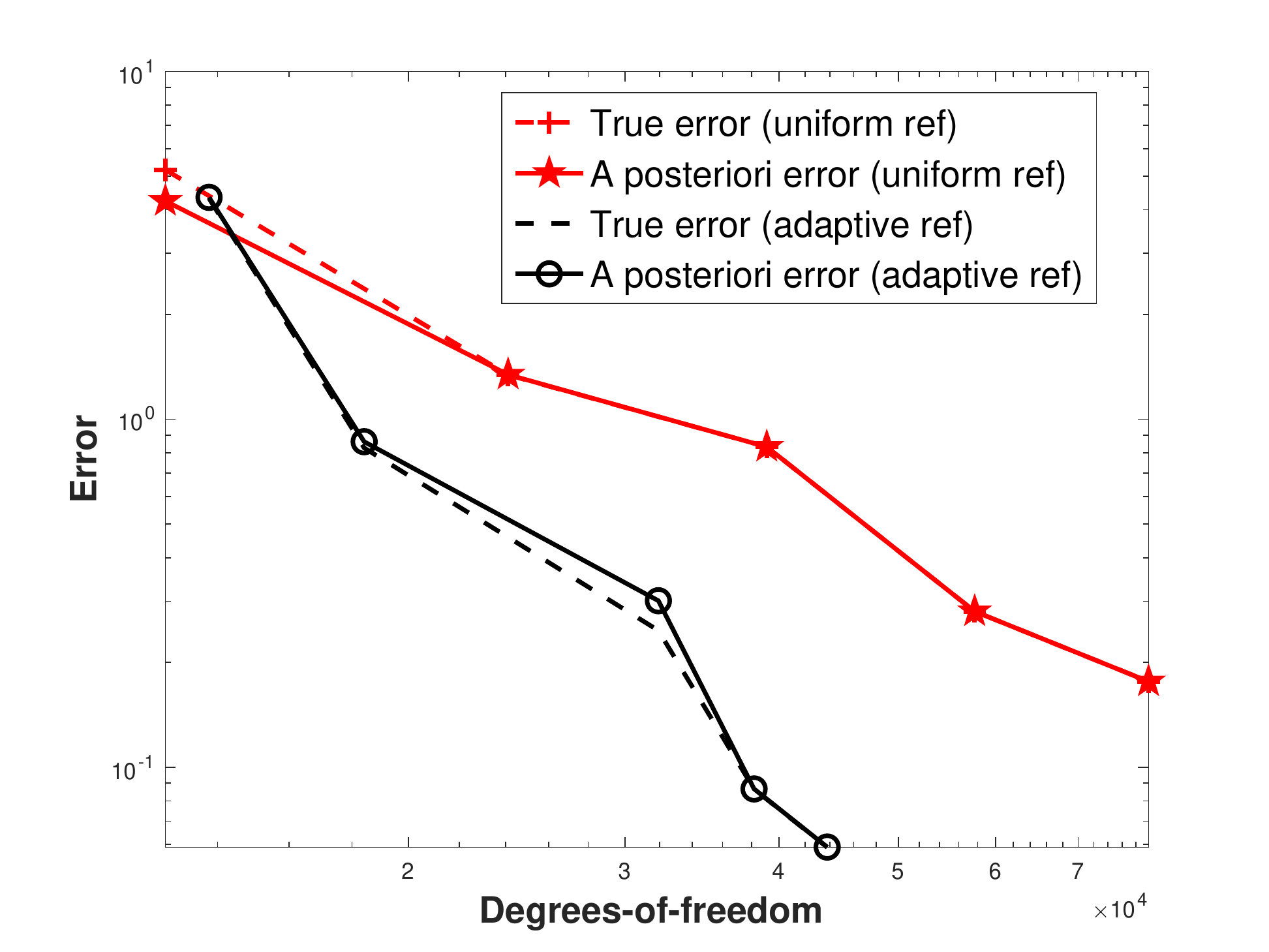}
\end{tabular}
\caption{Convergence comparison between uniform and adaptive p-refinements for the bump function solution.}
\label{fig:non_uniform_adaptive}
\end{center}
\end{figure}

\begin{figure}
\begin{center}
  \begin{tabular}{c}
  \includegraphics[angle=0,width=17.0cm,height=7cm]{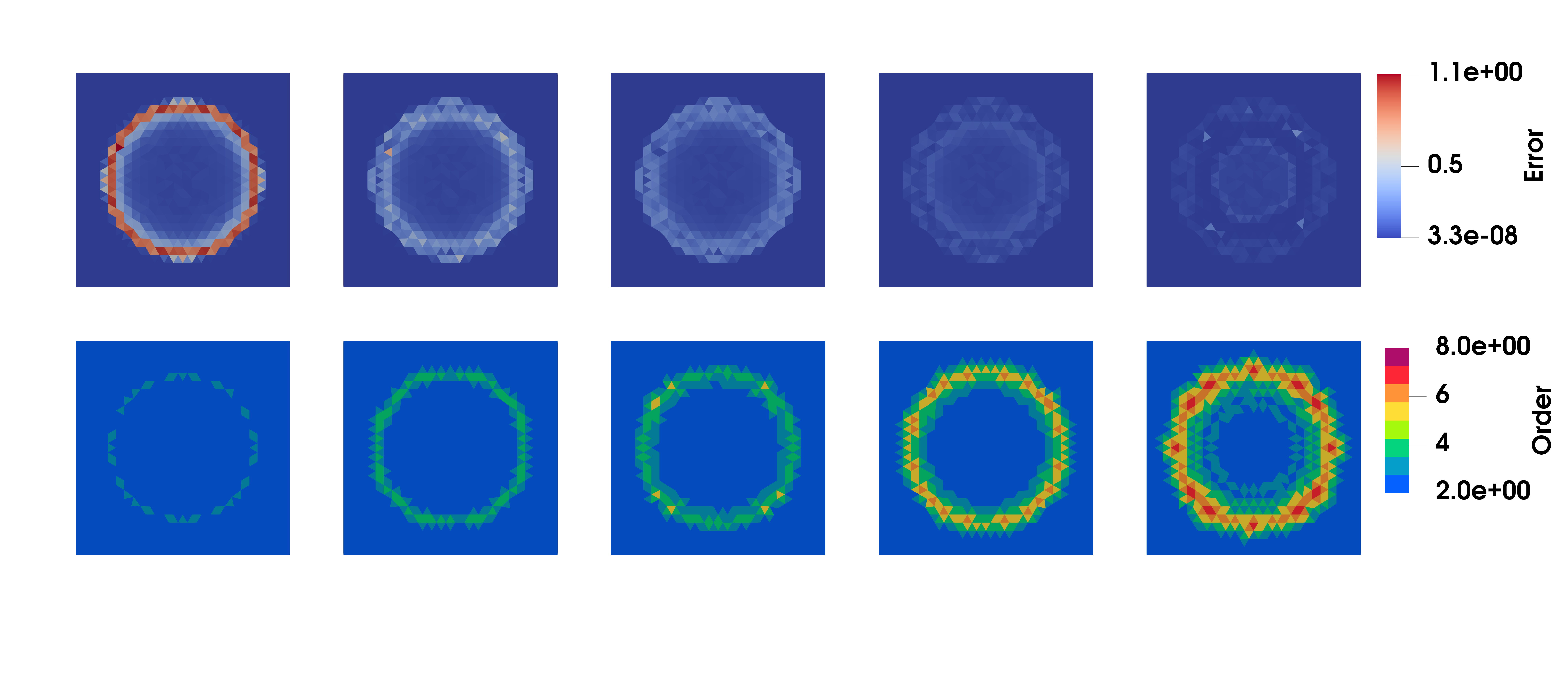}
\end{tabular}
\caption{Distribution of the error (top row) and polynomial order (bottom row) in a p-adaptive sequence.}
\label{fig:non_uniform_sequence}
\end{center}
\end{figure}

Thus it is expected that most of the error of the finite element approximation concentrates along
the circle. This becomes evident in the distributions of the error as well as polynomial order over the mesh as displayed
in the p-adaptive sequence as shown in Fig.~\ref{fig:non_uniform_sequence}. The corresponding convergence result, shown in
Fig. \ref{fig:non_uniform_adaptive}, exhibits the better performance compared to the uniform p-adaptive strategy.

\subsubsection*{Rough solution with singularities}\label{sec:rough-sol}
\begin{figure}
\begin{center}
 \begin{tabular}{cc}
 \includegraphics[angle=0,width=7.0cm,height=6.2cm ]{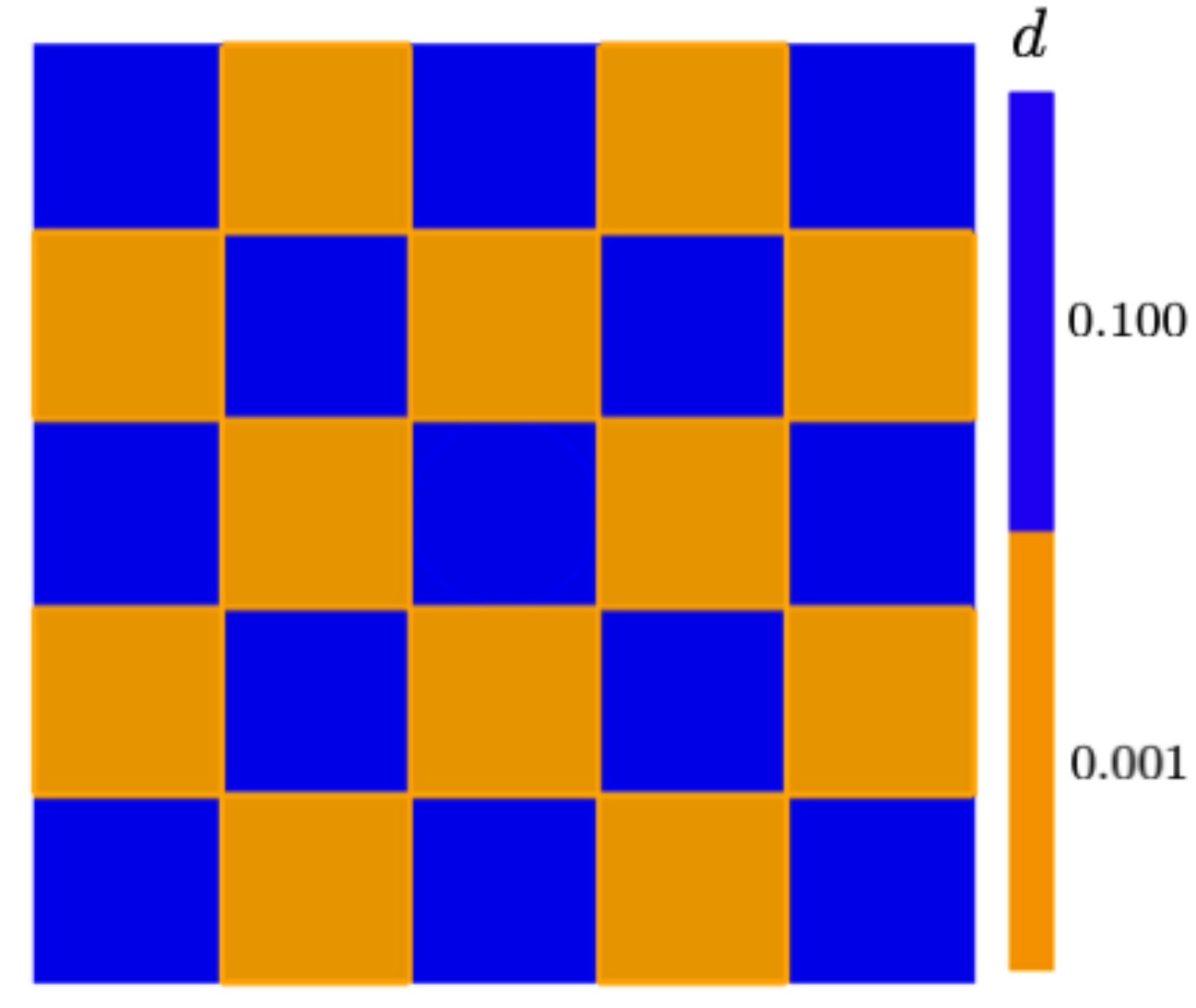} &
  \includegraphics[angle=0,width=7.0cm,height=6.3cm ]{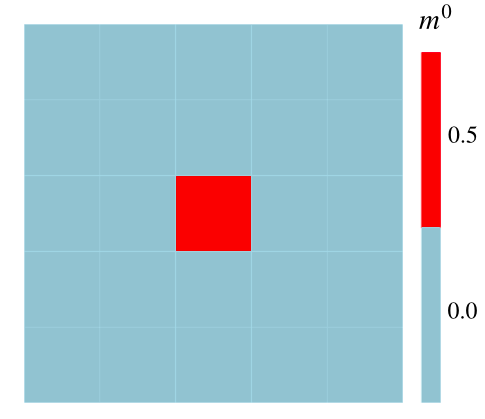}\\
  (a) & (b)
\end{tabular}
\caption{The heterogeneous domain $\Omega = [-1,~1]^2$ with a checkerboard pattern where in (a) the blue regions have
$d = 1\times10^{-3}$ and the rest have $d = 1\times10^{-1}$, and in (b) an initial value $m^0 = 0.5$ is prescribed.}
\label{fig:check-board}
\end{center}
\end{figure}

Consider a one-species reaction diffusion problem on the square domain $\Omega = [-1,~1]^2$, with reaction kinetics of
Fisher's type,
\begin{equation}\label{eq:fisher}
f = m[1-m].
\end{equation}
The domain is comprised of square patches with contrasting diffusivities $\bm{D} = d(\bm{x})\bm{I}$ with a checkerboard
pattern, as shown in Fig.~\ref{fig:check-board} (a). For the
blue patches $d = 0.1$ and for the remaining patches $d=0.001$. An initial condition of $m^0 = 0.5$ on the centre
square and zero elsewhere is prescribed - see Fig.~\ref{fig:check-board} (b). A homogeneous flux
boundary condition of the type  given in equation \eqref{eq:bd-E} is prescribed along the entire boundary $\Gamma$ (i.e., $\bar{h} = 0$).

\begin{figure}
\begin{center}
 \begin{tabular}{c}
 \includegraphics[angle=0,width=13.0cm,height=12.0cm ]{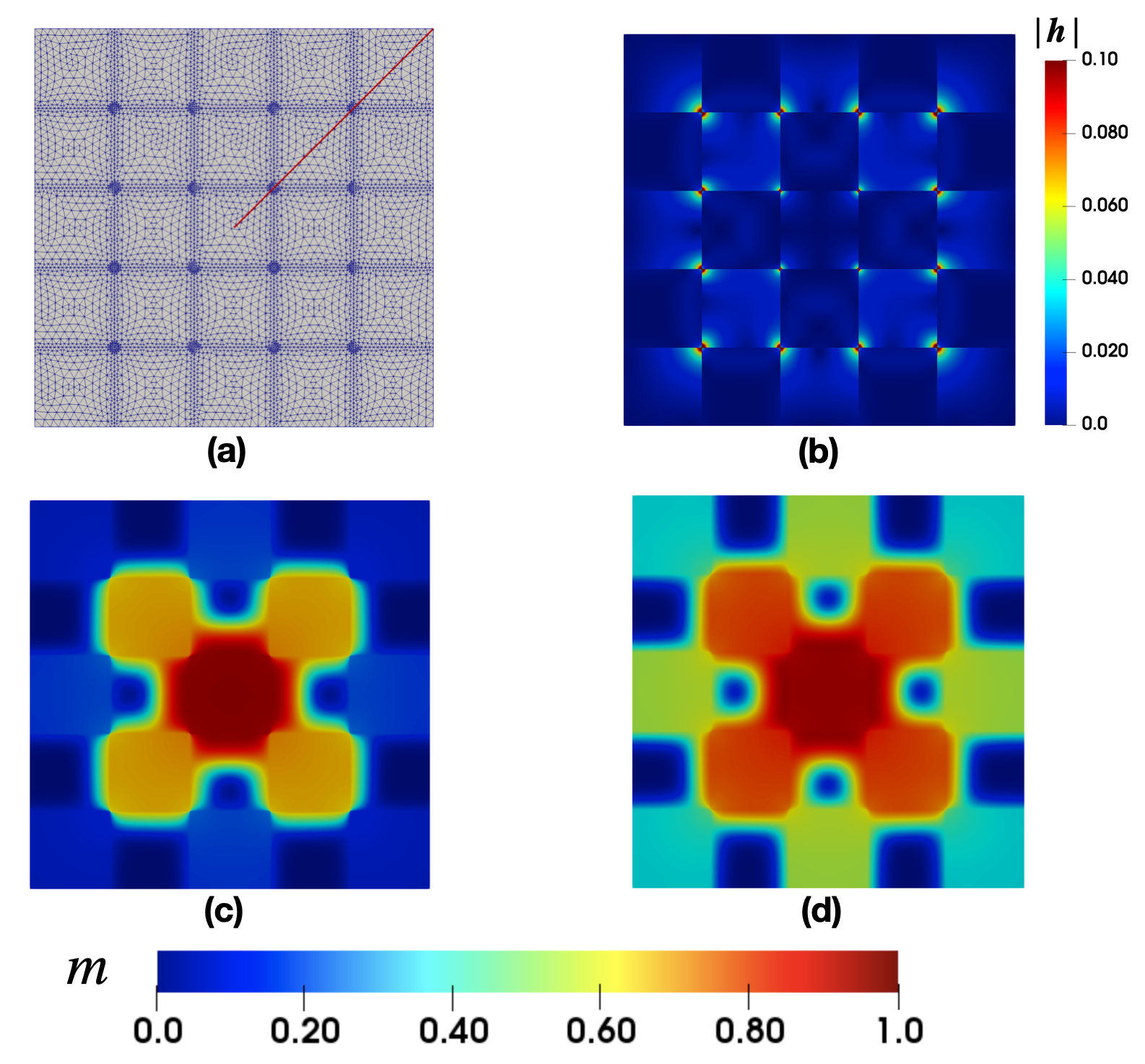}
\end{tabular}
\caption{Comparison of approximations by standard and mixed methods, (a) an a priori adaptively refined mesh on which the
mass concentration approximation at $t=6$ are computed based on the standard (c) and mixed (d) methods. (b) shows the distribution
of magnitude of flux, computed using the mixed method using the mesh (a) at $t=6$, where singularities at the
corners of the patches are shown.}
\label{fig:check-board-mass}
\end{center}
\end{figure}

Solutions were computed up to $t = 6$ with a time-step size $\Delta t_n = 0.1$. Both the standard and mixed solutions were
computed and compared. Because of the heterogeneous diffusivity, the solution develops kinks along the
interfaces of the patches and singularities at the corners. It is important to note the well-known fact that such irregularities
 (singularities) cannot be resolved by increasing the polynomial order. A feasible way of resolving such features is using
local h-adaptivity. In fact, numerical experimentation (not presented here) showed that p-adaptivity caused artificial
oscillation near the corners as the polynomial order increases locally. For example, Fig.~\ref{fig:check-board-mass} (b) shows
the distribution of the flux magnitude, computed using the mixed method with the a priori adaptively refined mesh as shown
in Fig.~\ref{fig:check-board-mass} (a). Figs.~\ref{fig:check-board-mass} (c) and (d) show the mass distribution
computed using the standard and mixed methods, respectively. Even though the reference mesh Fig.~\ref{fig:check-board-mass} (a)
has been used in both cases, the difference in their respective solutions is apparent. This is due to the fact that the
standard formulation uses a $H^1$-conforming space and is unable to approximate solutions with features such as
discontinuities and singularities. By contrast, the mixed formulation uses a non-conforming $L^2(\Omega)$ space for $m^h$ and
the exact, physically motivated conformity for the flux $\bm{h}^h$, i.e., $H(\Div)$. This allows the mixed method to capture
such features.

\begin{figure}
\begin{center}
 \begin{tabular}{c}
\includegraphics[angle=0,width=8.0cm,height=7.2cm ]{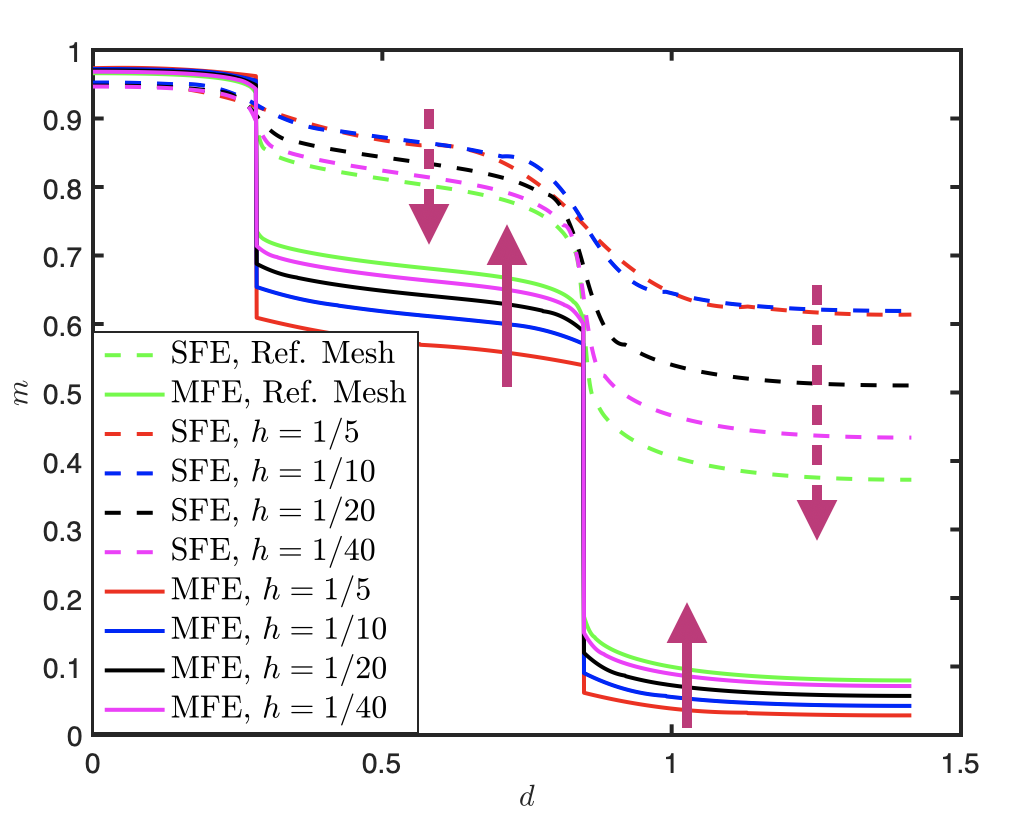}
\end{tabular}
\caption{Mass concentration profile along a diagonal line from the centre of the square to the right corner (as indicated
in Fig.~\ref{fig:check-board-mass}(a) with a red line segment) at $t=6$. Various uniform meshes at different refinement levels
described by the mesh parameter $h$ along with a reference mesh (designated as Ref. Mesh), which was adaptively refined
along the interfaces and corners of the patches, were used. }
\label{fig:check-board-plots}
\end{center}
\end{figure}

Fig.~\ref{fig:check-board-plots} shows the superiority of the
mixed method over the standard single-field formulation. Here, the mass  profiles along a line segment, coloured in red
in Fig.~\ref{fig:check-board-mass} (a), that connects the centre and the right top corner of the square domain $\Omega$
are displayed. Along this line segment, there are two interior corners of patches where discontinuities in $m_h$ are expected.
The solutions were obtained using various meshes at different levels of refinement with parameter $h = 1/5, 1/10, 1/20$, and
$1/40$, including a reference mesh (denoted by Ref. Mesh), which is obtained by refining along the interfaces and corners of
the patches. The discontinuities at the interior corners are captured almost exactly using the mixed method
regardless of the refinement level. By contrast, using the standard formulation, none of the meshes resulted in a
reasonable approximation of the discontinuities. Another interesting observation is that the approximation of the
mixed formulation converges from below. This is opposite to the standard formulation in which the
approximation overestimates the solution.

\subsection{Speed of travelling wave solutions}\label{sec:wave-speed}
The one-species Fisher-type equation \eqref{eq:fisher} supports travelling wave solutions. The wave nature of the solution depends on
the relative size of the reaction and diffusion terms. When the reaction term is dominant, the
wave front steepens and the wave travels with a finite speed. By contrast, when diffusion is dominant, the influence of the
reaction term becomes less and the solution exhibits typical diffusion behaviour, that is, it decays exponentially.

\begin{figure}
\begin{center}
\begin{tabular}{c}
\includegraphics[angle=0,width=11.0cm,height=5.4cm ]{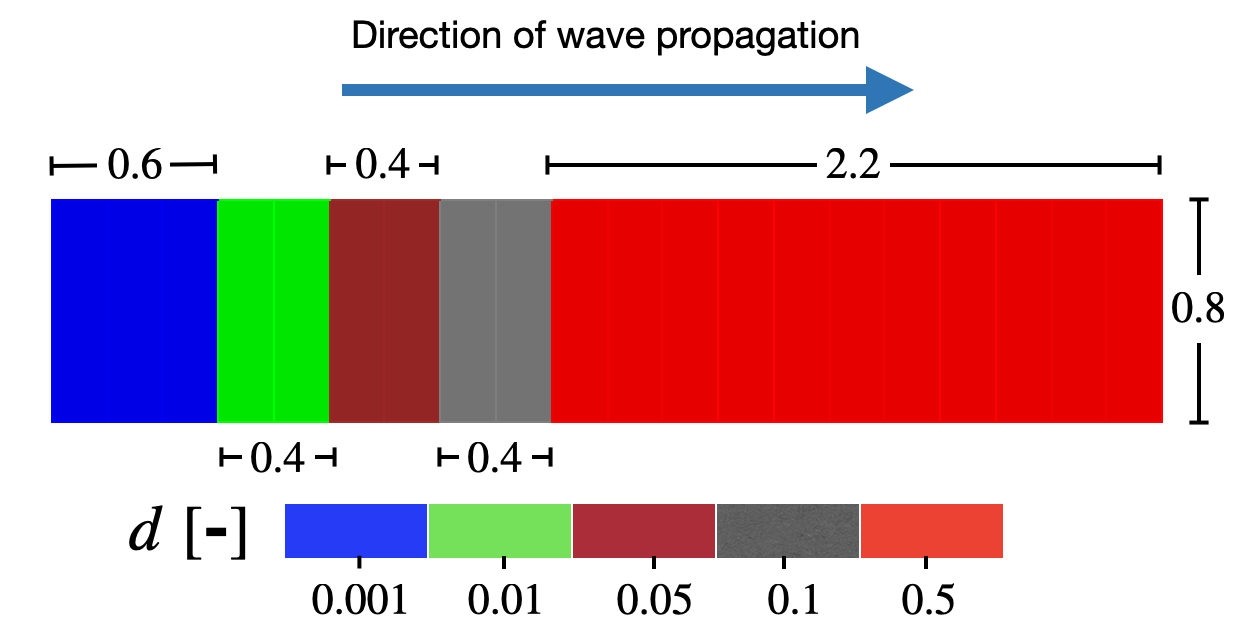}
\end{tabular}
\caption{Domain for the travelling wave problem. Color indicates the diffusivity, which is piecewise constant and increasing towards the right.}
\label{fig:rect_patch}
\end{center}
\end{figure}

\begin{figure}
\begin{center}
\begin{tabular}{c}
\includegraphics[angle=0,width=12.0cm,height=10.cm ]{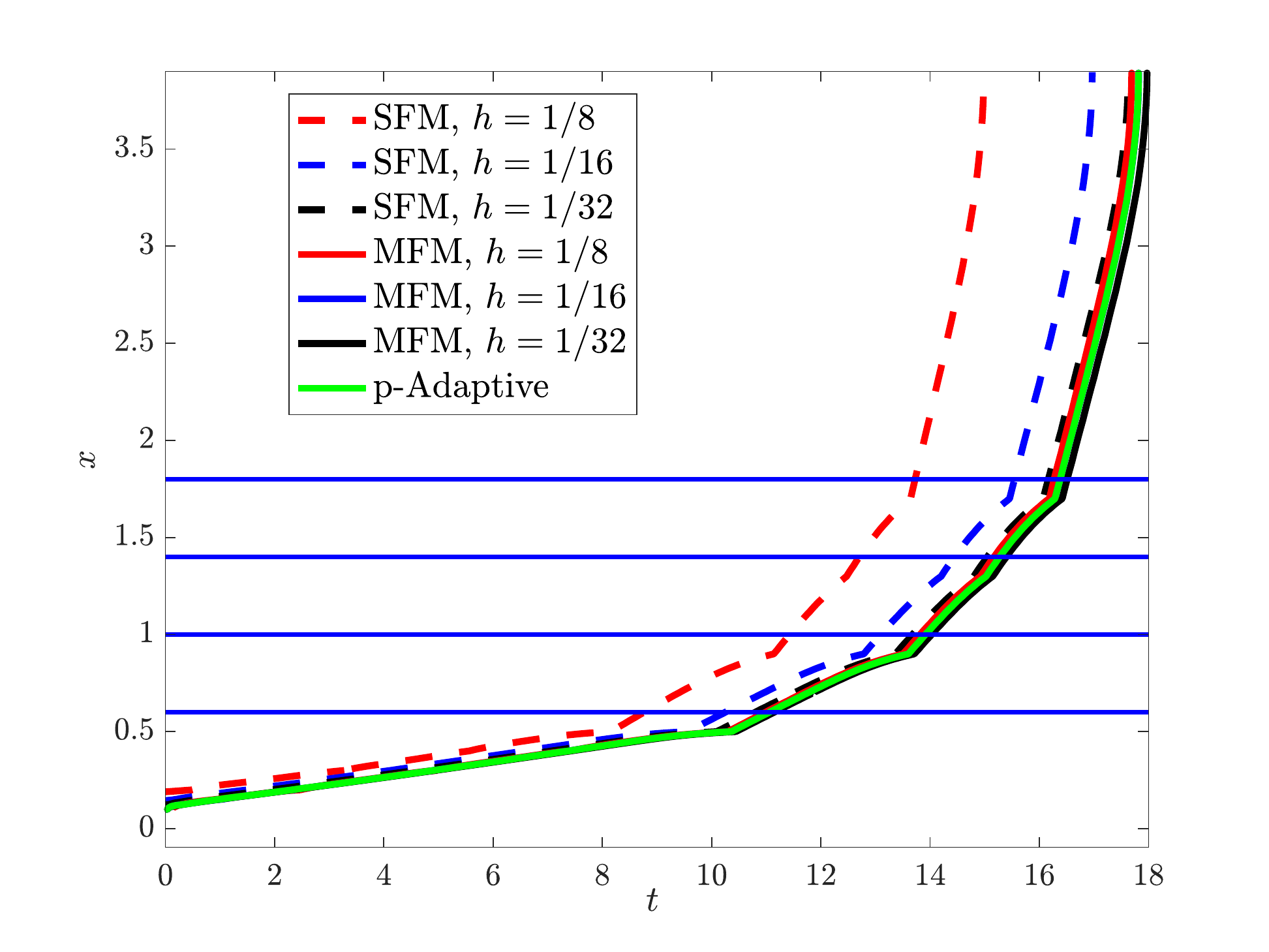}
\end{tabular}
\caption{Comparison of speed of travelling wave solution with time. The broken lines correspond to the computation
using the standard formulation (SFM) while the solid lines correspond to the mixed formulation (MFM).}
\label{fig:speed_comparison}
\end{center}
\end{figure}

Consider the planar domain as shown in Fig.~\ref{fig:rect_patch}, composed of rectangular patches that are arranged horizontally
with increasing diffusivity between successive patches. The rectangular domain has height $0.8$ and the width of the 5
patches are $0.6$, $0.4$, $0.4$, $0.4$, and $2.2$. The maximum diffusivity $d = 0.5$ is chosen so that the problem remains in the wave
propagation regime. As the initial condition, the mass concentration $m = 1$ is set on the first one-third of the left-most patch,
while on the rest of the domain $m$ is set to zero at $t=0$. Since the diffusivity is smallest in the first left patch, the solution
starts to evolve slowly with a sharp wave front. As the wave passes each interface its speed increases while the sharpness of the wavefront
decreases. The wavefront is identified using a levelset method based on a mass concentration value of $0.6$. Uniform time steps of 
length $\Delta t = 0.05$ have been used. The position of the wavefront
against time is presented in Fig.~\ref{fig:speed_comparison} for
various levels of mesh refinement. Importantly, the approximations based on the mixed method are accurate and converge to the correct
wave speed. However, the approximation using the standard method initially overestimates the true
speed of the travelling wave solution, and only slowly converges to the correct speed. The p-adaptive mixed formulation is used with the second
coarsest mesh, and it is also obtained that the speed of the wave is in good agreement with that of the converged results of the mixed or 
standard methods.  

\subsection{Pattern formation in ecological applications}
The class of problems considered here have application in various important areas including  biological pattern formation,
morphogenesis \cite{GARIKIPATI2017192} and electrophysiology \cite{Goktepe2009}.

\subsubsection*{Segregation pattern} A competition-diffusion model involving three interacting species is considered. The level and mode
of interaction between the species is the same. This, in effect, means that the magnitude of each species that is consumed by the others
is the same as the other species that consumes it. The reaction term, for each $i = 1, 2, 3$, is given by
\begin{equation}\label{eq:three-species}
f_i = \,m_i[1-a_{i1}m_1 - a_{i2}m_2 - a_{i3}m_3],
\end{equation}
\begin{figure}
\begin{center}
\begin{tabular}{c}
\includegraphics[angle=0,width=17.6cm,height=11.2cm ]{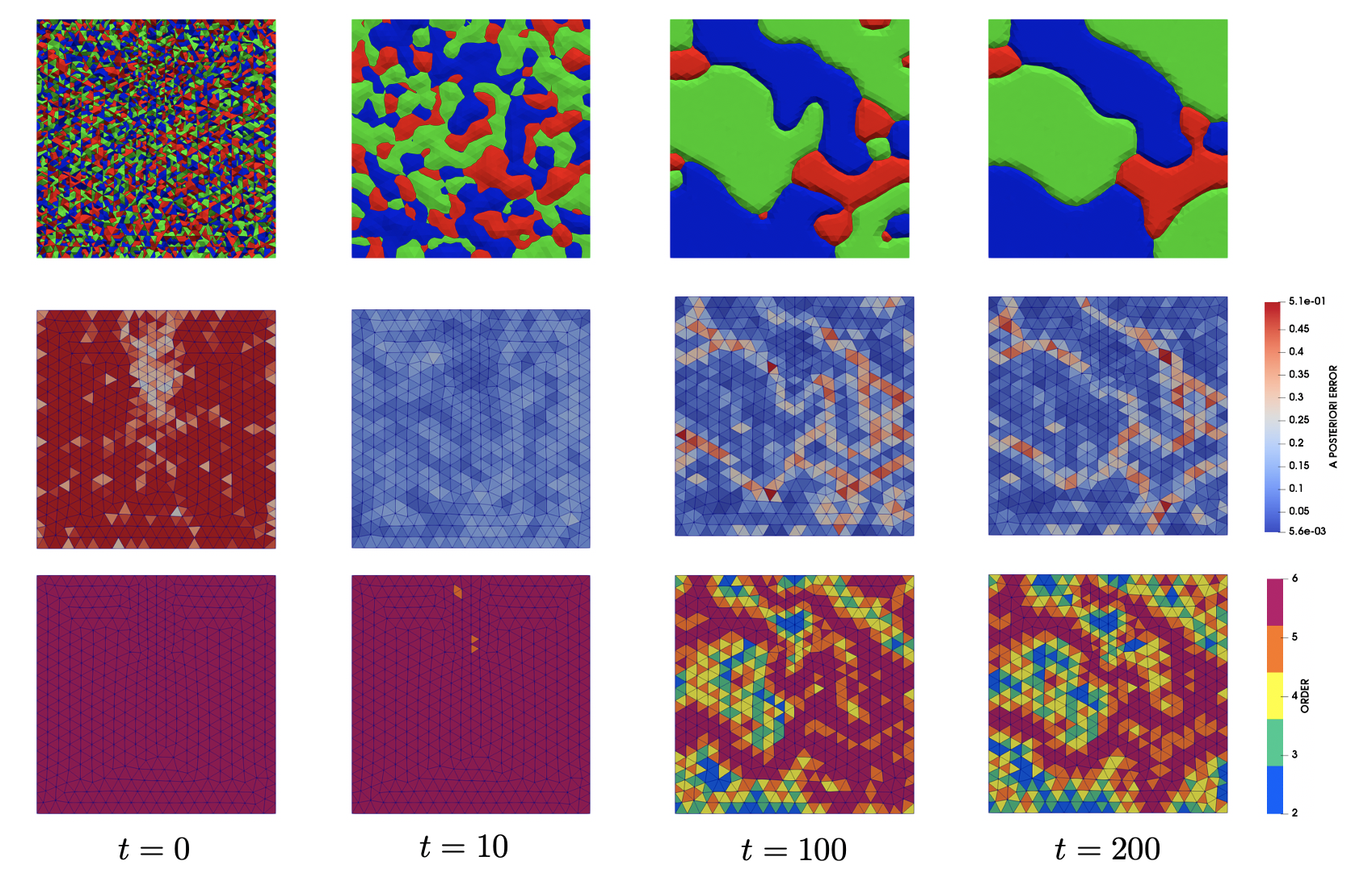}
\end{tabular}
\caption{Development of segregation pattern as a result of interactions of the three species (Equation \eqref{eq:three-species}, Table~\ref{tab:param-3})
at various times. The blue color surface plot represents the region
dominated by $m_1$ (i.e., $m_1>m_2$ and $m_1>m_3$), red color surface plot by $m_2$, and yellow surface plot by $m_3$.}
\label{fig:segregation}
\end{center}
\end{figure}
where the parameters in the model are represented in Table~\ref{tab:segregation}. It is assumed that all the three species have
the same mobility rates, i.e., $\bm{D} = d\bm{I}$, where $d = 0.01$. In cases where the dynamics is largely influenced by the
reaction term, it is important to analyse the local stability of the spatially homogeneous problem (i.e., ignoring the diffusion terms).
Such an analysis provides important insights into the range of parameter values for various possible spatio-temporal interaction patterns.
A local stability analysis of the problem described by equations \eqref{eq:three-species} reveals eight equilibrium points of which
only, namely $[m_1,~m_2,~m_3]_1=[1/a_{11},~0,~0]$, and $[m_1,~m_2,~m_3]_2=[0,~1/a_{22},~0]$, and
$[m_1,~m_2,~m_3]_3=[0,~0,~1/a_{33}]$ are quasi-stable. In the case of the segregation problem, these equilibrium points represent
regions which are exclusively occupied by one of the species.
\begin{table}
\begin{center}
\caption{List of parameters for three species segregation problem} \label{tab:param-3}
\begin{tabular}{ccc|ccc|ccc|ccc}
\hline
$d_1$ & $d_2$ & $d_3$ & $a_{11}$ & $a_{12}$ & $a_{13}$ & $a_{21}$ & $a_{22}$ & $a_{23}$ & $a_{31}$ & $a_{32}$ & $a_{33}$ \\
\hline
$0.01$ & $0.01$ & $0.01$ & $1$ & $3$ & $3$ & $3$ & $1$ & $3$ & $3$ & $3$ & $1$
\end{tabular}
\label{tab:segregation}
\end{center}
\end{table}

It is assumed initially that all the three species are distributed randomly over the domain $\Omega = [-1,~1]^2$ as shown in
Fig.~\ref{fig:segregation}. Such problems have been studied previously \cite{Mergia2020, Mimura2015, Mimura1986}.
In those studies the numerical approaches were either the finite  difference method or the standard finite element method.

A relatively coarse mesh for such a problem is used, but to sufficiently represent the random initial condition and to captured the fast
dynamics at the beginning stage order 6 polynomial approximation are used as shown in the lower left corner of Fig. \ref{fig:segregation}.
During the early stages of evolution, as show in Fig.~\ref{fig:segregation} at $t=0$ through $t=10$, the dynamics appears to be 
reasonably fast. Eventually, as shown in the second row of Fig.~\ref{fig:segregation}, as the species start to establish themselves 
into well defined regions each occupied by one of the species, the interaction starts to proceed in a slower manner. During this time
the error distribution becomes more concentrated in the vicinity of the boundaries of these regions. It is shown that the polynomial 
adaptation follows the error distribution very closely. These regions tend to a quasi-stable configuration, that is patches of convex 
shapes with triple junctions with angle of separation given by $2\pi/3$.
\subsubsection*{Cyclic interaction}
Here a three-species competition-diffusion system of equations with the reaction term given by  equation \eqref{eq:three-species}
is considered. The species react with each other in a cyclic way (based on parameters in Table~\ref{tab:param-cyclic}) resulting in
various complex spatio-temporal patterns such as spiral-like, and band-like structures depending on the topology of the habitat
and the initial configuration.
\begin{table}
\begin{center}
\caption{List of parameters for three species cyclic interaction problem} \label{tab:param-cyclic}
 \begin{tabular}{ccc|ccc|ccc|ccc}
\hline
$d_1$ & $d_2$ & $d_3$ & $a_{11}$ & $a_{12}$ & $a_{13}$ & $a_{21}$ & $a_{22}$ & $a_{23}$ & $a_{31}$ & $a_{32}$ & $a_{33}$ \\
\hline
$0.01$ & $0.01$ & $0.01$ & $1$ & $2$ & $7$ & $7$ & $1$ & $2$ & $2$ & $7$ & $1$
\end{tabular}
\end{center}
\end{table}

\begin{figure}
\begin{center}
 \begin{tabular}{c}
 \includegraphics[angle=0,width=16.6cm,height=11.2cm ]{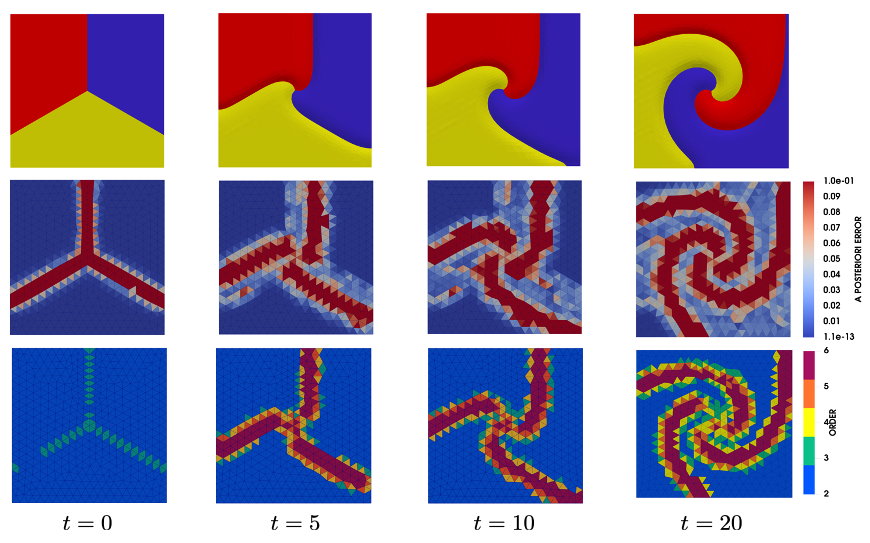}
\end{tabular}
\caption{Development of spiral pattern as a result of cyclic interactions of three species (Equation \eqref{eq:three-species}, Table~\ref{tab:param-cyclic})
at various stages. Region with blue depicts $m_1$, red $m_2$, and green $m_3$.}
\label{fig:three-species}
\end{center}
\end{figure}

A square habitat $\Omega = [-1,~1]^2$ is considered with initial configuration as shown in Fig.~\ref{fig:three-species} at  $t=0$ (top row left).
The parameter values considered in the simulations are presented in Table~\ref{tab:param-cyclic}. A relatively coarse mesh is used for such
problem whose solution have complex, fine and spiral structure. A uniform time stepping is also used with step size $\Delta t = 0.2$. Initially the
polynomial order is set to 2 and adaptively increases to 6 (with parameter $\theta_\mathrm{max} = 0.6$ and $\theta_\mathrm{min} = 0.1$) 
throughout the simulation. As shown in Fig. \ref{fig:three-species}, a spiral pattern starts to form turning in a clockwise direction. The 
spiral shape consists of stripes of each species lying side-by-side, and eventually fills the region and continues with the same spiral feature. 
It is also shown that
the error is high in the vicinity of the interface between the species which led to the adaptivity taking place only on elements around 
such interfaces. This clearly shows the efficiency of the p-adaptive mixed method.  

\subsection{Spiral wave re-entry in electrophysiology}

The propagation of ionic current in the cardiac muscle can be simulated using a monodomain model, which is mathematically equivalent to
the reaction-diffusion equation. The transmembrane electric potential can be viewed as a diffusing species, which ``reacts''
locally with the cellular ion channel densities. The reaction term depends on the ion channel densities through a set of highly
nonlinear and coupled ordinary differential equations. Thus, the ion channel densities can be viewed as non-diffusing species and
treated simply as internal state variables. There are a large number of models available for the reaction term (called cardiac
electrophysiology models) with varying degrees of complexity in terms of the number of (internal) variables. Using the proposed
mixed method, we simulate the phenomenon of spiral wave re-entry --- the cause of several cardiac arrhythmias, such as ventricular
tachycardia, atrial flutter, and atrial and ventricular fibrillation.

A square block of cardiac tissue of dimension $100\;\mathrm{mm}$ is considered. The domain is subdivided into a relatively coarse 
triangular mesh. The propagation of the transmembrane electric potential (more commonly known as the action potential) is governed by
\begin{equation}\label{eq:aliev-panfilov}
\frac{\mathrm{d}}{\mathrm{d}\tau}m - \mathrm{div}(\bm{D}\nabla m) =  f(m,r) \,+ \,I_{\Omega'}(\tau),\\
\end{equation}
where the non-dimensional variable $m\in[0,1]$ is related to the transmembrane action potential $E~[\mathrm{mV}]$ through the relation
\[
E = [100 m - 80]~\mathrm{mV},
\]
$r$ is a single internal variable representing the density of ionic channels, and $I_{\Omega'}$ is the external stimulus. The
time $t~[\mathrm{ms}]$ is non-dimensionalised as
\[
t = 12.9 \tau~\mathrm{ms}.
\]
One of the simplest models capable of reproducing the spiral wave re-entry, the Aliev-Panfilov model \cite{Rubin1996}, is adopted
for the reaction term:
\begin{equation}\label{eq:aliev-panfilov2}
f(m,r) = cm[m-\alpha][1-m] - rm.
\end{equation}
Equation \eqref{eq:aliev-panfilov2} is supplemented by an ordinary differential equation for the recovery (internal state) variables $r$:
\begin{equation}\label{eq:recovery}
  \frac{\mathrm{d}}{\mathrm{d}\tau}r = \bigg[\gamma + \frac{\mu_1r}{\mu_2+m}\bigg][-r-cm[m-b-1]],
  \end{equation}
where the parameters appearing in equations \eqref{eq:aliev-panfilov2} and \eqref{eq:recovery} are given in Table~\ref{tab:param-reentry}.
We assume the conductivity to be isotropic, i.e., $\bm{D} = d\bm{I}$. The simulation is carried out using the IMEX mixed formulation
with order $k=1$ (recall that $k$ is the polynomial order used for the approximation of $m$ and $k+1$ is the order used for the flux). The 
p-adaptivity strategy with parameters $\theta_\mathrm{max} = 0.7$ and $\theta_\mathrm{min} = 0.03$ starts with uniformly order 2 and increases
locally to order 6.

A horizontal planar wave is initiated by setting the action potential to $E = -40\;\mathrm{mV}$ on the region between $y = 0$ and $y = 3$. 
The wave form continues to propagate upwards as seen from the snapshot at $t=160\;\mathrm{ms}$. 

\begin{figure}
\begin{center}
\begin{tabular}{c}
\includegraphics[angle=0,width=16.6cm,height=8.2cm ]{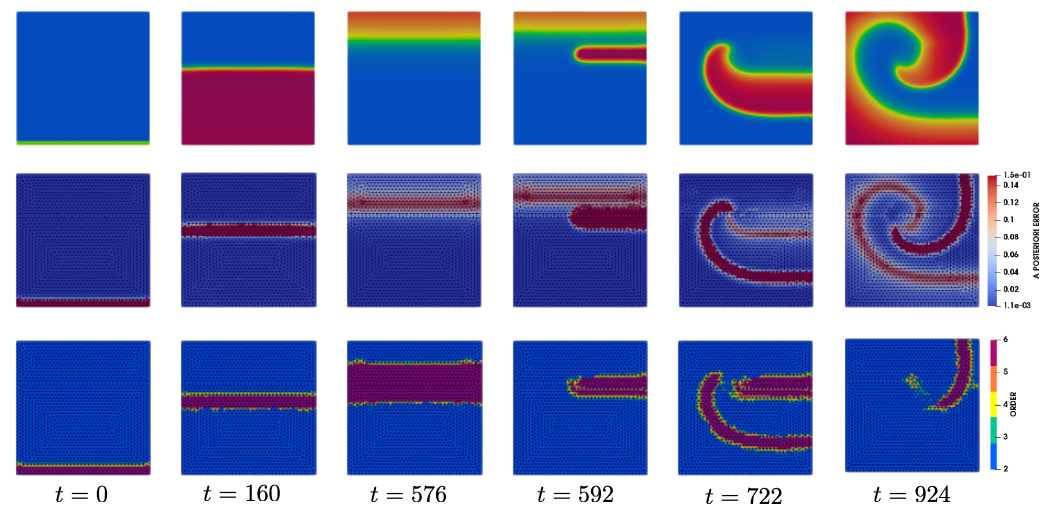}
\end{tabular}
\caption{Evolution of a planar wave into a rotating spiral wave re-entry as a result of external stimulation $I_{\Omega'} = 40$ applied on the
red shaded region (top left) during the time interval $t = 570~\mathrm{ms}$ to $t = 580~\mathrm{ms}$. The planar wave is initiated with
an initial excitation of $E(t=0) = -40~\mathrm{mV}$ on the region shaded in blue.}

\label{fig:wave-reentry}
\end{center}
\end{figure}

Before the depolarising tail disappears, an external stimulus $I_{\Omega'}$ is applied to the strip of region, defined by 
$\Omega' = \{(x, y): 50 < x < 100,~67 < y < 70\}$, in order to initiate the spiral wave re-entry. The stimulus has magnitude $40$ 
and is applied at $t=565~\mathrm{ms}$ for a duration of $10\;\mathrm{ms}$. This results in the development of the wavebreak 
(shown in Fig.~\ref{fig:wave-reentry} at $t = 7722~\mathrm{ms}$). The wavebreak then evolves into a stable rotating vortex, as shown in 
the snapshots at $t=924\;\mathrm{ms}$, and continues afterwards.

\begin{table}
\begin{center}
\caption{List of parameters for spiral wave re-entry problem} \label{tab:param-reentry}
\begin{tabular}{ccccccc}
\hline
$d$ & $\alpha$ & $\gamma$ & $b$ & $c$ & $\mu$ & $\mu_2$ \\
 $[\mathrm{mm}^2]$ & $[-]$ &$[-]$&$[-]$&$[-]$&$[-]$&$[-]$\\
\hline
$0.01$ & $0.01$ & $0.01$ & $1$ & $2$ & $7$ & $7$
\end{tabular}
\end{center}
\end{table}
The computational aspects of this problem have been considered by several researchers in the electrophysiology and
electromechanics community. One notable work is by G\"oktepe and Kuhl \cite{Goktepe2009} in
which they used the standard finite element approach with an implicit time integration scheme on a structured quad mesh. The proposed p-adaptive, 
IEMX mixed, as compared to their work, 
method is seen to capture the spiral wave re-entry dynamics more efficiently.

\section{Conclusion}
A family of p-adaptive implicit-explicit, mixed finite element formulations has been proposed for a general class of reaction-diffusion based problems.
In contrast to single-field, standard finite element formulations, this class of methods provides accurate approximations of a wider
class of solutions, including those with less regularity. A standard formulation was shown to converge poorly, if at all, for such problems.
The IMEX approach has been shown to be efficient and
eliminates the dependence of algorithmic stability on the size of the spatial mesh size by handling the non-local diffusion part implicitly.
This  advantageous feature allows for mesh refinement, for example in an adaptive strategy, without the need for changing the time step size,
$\Delta t$. The explicit treatment of the local reaction term makes the implementation generic and modular for various classes of reaction
kinetics as demonstrated by the wide range of problems that have been analysed in this paper.
The finite element spaces are built using a hierarchical construction which, in addition to offering optimal conditioning of the resulting
linear system, makes the use of the $p$-adaptivity strategy a natural choice \cite{ainsworth2003hierarchic}. The mixed formulation introduces
additional DoFs. However, the computational complexity due to this increase in DoFs can be handled efficiently using
static condensation as the mass concentration field (which is in $L^2$)
can be inverted locally since the local contributions are decoupled from one another. Moreover, this local inversion can also be
used in block iterative schemes that involves computation of the Schur complement as an intermediate step. The Schur complement can
then be computed exactly, resulting in a sparse global structure, rather than reverting to the common practice of approximating it.

A distinguishing feature of the mixed formulation is that it leads to straightforward derivation and implemetation of  residual based a
posteriori error estimations without the need for computationally demanding postprocessing effort as it is usually the case in 
literature, see, for example, \cite{Ainsworth2008, Brass1996, Lloyd2013}. This feature together with the hierarchical approximation of 
the mixed finite element method is exploited in formulating the p-adaptive strategy. It has been demonstrated 
by a range of examples that the p-adaptive algorithm performs very well in efficiently resolving fine features.  

The performance of the proposed
formulation is demonstrated by a number of challenging examples. The advantages of the proposed method over the standard techniques
are showcased by the following two examples: a problem that has singularities (see Section~\ref{sec:rough-sol});
and one that supports travelling wave solutions (see Section~\ref{sec:wave-speed}). The capability of this general p-adaptive framework is
demonstrated by applying it to problems arising from different applications such as electrophysiology \cite{Goktepe2009, Rubin1996},
and spatial pattern formation in theoretical ecology \cite{Mimura1986, Morishita2008, Mimura2015}.

Through a generalisation straightforwardly the proposed mixed method can be coupled
to the mechanical deformation field for applications in cardiac electromechanics and chemo-mechanics. Due to the explicit treatment of the
reaction term, the approach can be easily linked with, for example, electrophysiology models in the CellML repository \cite{Lloyd2013}.
Similarly, due to this explicit treatment of the reaction term, our computational approach can be easily used to drive the form of reaction
kinetics models from experimental data following the approach proposed in \cite{Brunton3932}.
\section*{Acknowledgements}
The authors gratefully acknowledge the support provided by the EPSRC Strategic Support Package: Engineering of Active Materials
by Multiscale/Multiphysics Computational Mechanics - EP/R008531/1.
\section*{References}

\bibliography{mybibfile1.bib}

\end{document}